\documentclass[preprint,showkeys]{revtex4-1}
\usepackage{graphicx,subfigure}
\usepackage{amsmath,amsfonts,amssymb,amsbsy,amscd,amsgen,amsthm}
\usepackage[]{natbib}
\usepackage{color}
\usepackage{appendix}
\usepackage[colorlinks,citecolor=blue]{hyperref}
\usepackage{comment}

\newcommand{\id}{\mathrm d}
\newcommand{\vc}{\mathbf}
\newcommand{\bnabla}{\pmb\nabla}

\newcommand{\lins}[1]{\mathbf L_{#1}}

\newcommand{\etab}{\pmb\eta}

\renewcommand{\i}{\mbox{i}}

\theoremstyle{definition}

\newtheorem{defn}{Definition}

\begin{document}
\title{Dynamical indicators for the prediction of bursting phenomena in high-dimensional systems}
\author{Mohammad Farazmand} 
\author{Themistoklis P. Sapsis}
\affiliation{Department of Mechanical Engineering,
Massachusetts Institute of Technology,
77 Massachusetts Ave., Cambridge, MA, USA}

\begin{abstract}
Drawing upon the bursting mechanism in slow-fast systems, we propose  indicators for the 
prediction of such rare extreme events which do not require a priori known slow and fast 
coordinates. The indicators are associated with  functionals defined in terms of Optimally Time Dependent (OTD) modes. One such functional has the form
of the largest eigenvalue of the symmetric part of the
linearized dynamics reduced to these modes.
In contrast to other choices of subspaces, the proposed modes are flow invariant and therefore
a projection onto them is dynamically meaningful. We illustrate the application of these indicators on three examples: a prototype low-dimensional model, a body forced turbulent fluid flow, and
a unidirectional model of nonlinear water waves. We 
use Bayesian statistics to quantify the predictive power of the proposed indicators.
\end{abstract}
\keywords{rare events; probabilistic prediction; Kolmogorov flow; modified nonlinear Schr\"odinger equation; intermittency.}
\maketitle

\section{Introduction}
Complex irregular behavior is a characteristic of chaotic systems, which is
usually visualized through the time series of an observable. Many natural and engineering systems exhibit 
a second level of complexity typified by rare extreme bursts in the time series of 
certain observables. They are rare in the sense that they are short-lived
and the frequency at which they occur is significantly smaller than
the typical frequency of the time series; and
they are extreme in the sense that they correspond to the values of the observable
that are several standard deviations away from its mean value.
Examples of such rare, extreme phenomena in nature include
oceanic rogue waves~\citep{dysthe08, muller}, 
intermittent fluctuations in turbulent models~\citep{frisch80, Majda2014a, cousins_sapsis}
and extreme events in climate dynamics~\citep{cai14, chen_majda_giannakis}.
While the prediction of extreme events is of utmost importance, 
our dim understanding of their origins and precursors 
has impeded our ability to predict them. 

A promising approach is to predict the rare events directly from the time series of the
observable. If the system has a compact, finite-dimensional attractor, 
the dynamics can in principle be reconstructed from the observations by 
delay-coordinate embedding techniques~\citep{takens,sauer91, Berry2015}, or linear and/or nonlinear order reduction methods ~\cite{Williams2013, Chiavazzo2014, zhao15, Sargsyan2015, Zhong16}. 
However, for high-dimensional chaotic attractors the reconstructed dynamics have a poor forecasting skill (see e.g. ~\cite{Zhong16, Berry2016}) which is comparable with Mean Square Models (models based on carefully tuned Langevin equations \cite{Majda_filter}).  Since rare extreme events are associated with strong nonlinearities and intermittently positive Lyapunov exponents (i.e., high sensitivity to perturbations), 
their prediction from a finite set of observations is challenging and
remains an active area of research (see, e.g.,~\citet{majda12,BAK15}).

\begin{figure}
\centering
\includegraphics[width=.7\textwidth]{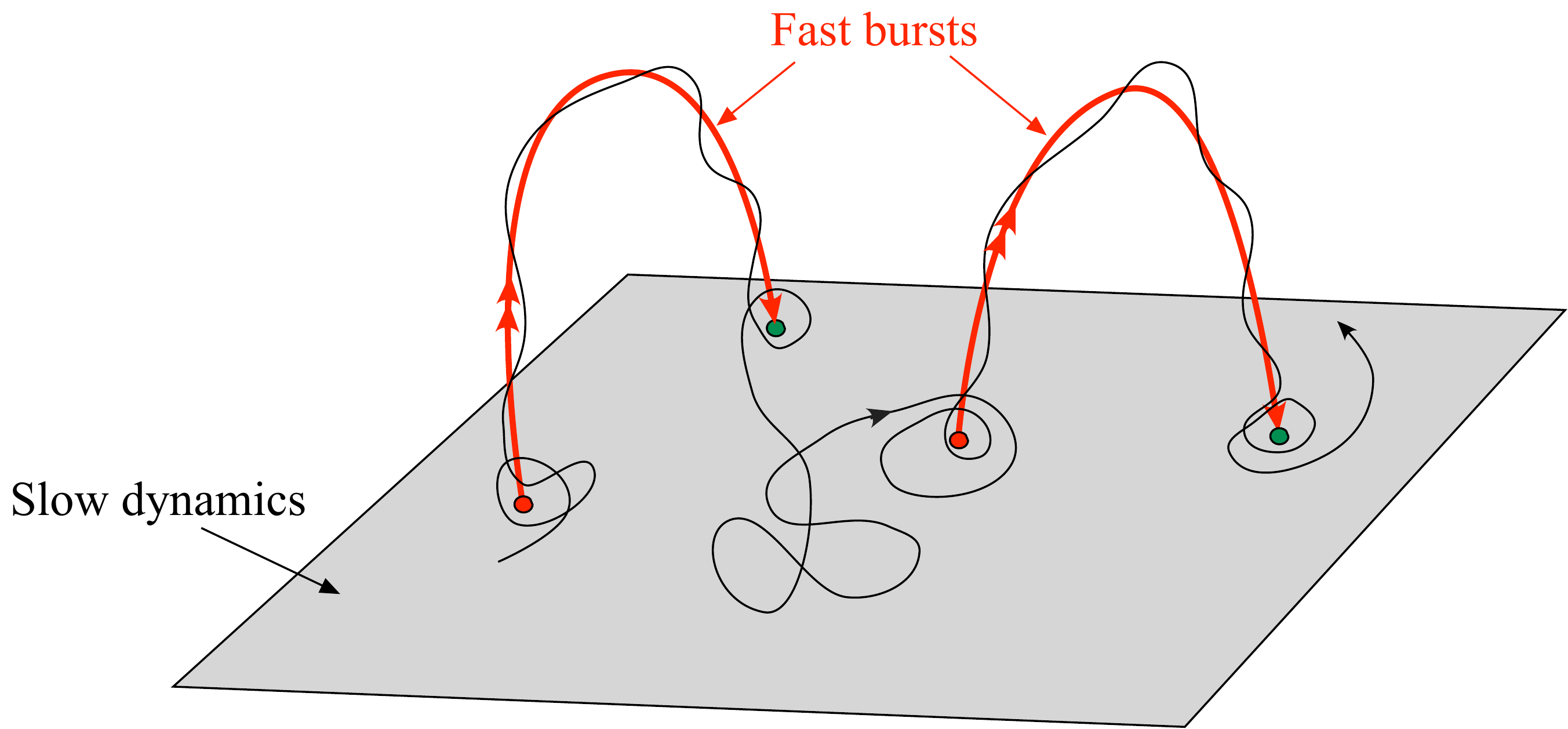}
\caption{An illustration of slow-fast systems with bursting orbits homoclinic
to the slow manifold. While we depict the slow manifold with a plane, it can in reality be 
a complicated high dimensional manifold.}
\label{fig:schem_bursting}
\end{figure}

A more physically illuminating approach comes from multiscale analysis, where
a dynamical system model is decomposed
into slow and fast variables~\citep{fenichel1979geometric,guckenheimer,wiggins1994normally} or stable and unstable manifolds ~\citep{Kevrekidis1987, Kevrekidis90}. The bursting mechanism in these models is rather 
well-understood~\citep{guck12,guck15}. For the most part, the dynamics takes place on the 
slow manifold. The slow dynamics may be chaotic, but no bursting 
event occurs on the slow manifold itself. The bursts occur along 
the unstable manifold (of the slow manifold) and correspond to the growth of the
fast variables. The unstable manifold is typically homoclinic to the
slow manifold such that the flow returns eventually to the slow manifold~\cite{haller99}.
This cycle can repeat indefinitely and, if the slow dynamics is
chaotic, irregularly (see figure~\ref{fig:schem_bursting}, for an illustration).

While this geometric approach is certainly illuminating, 
it is of little applicability to complex systems, since
a clear separation of time scales is often not available
in realistic models (e.g., Navier--Stokes equations). 
Nor does there exist a general recipe to transform the coordinates into slow and fast 
variables~\citep{pope06}. This becomes particularly prohibitive in high dimensional systems.

Here, we introduce a diagnostic indicator for the prediction of rare extreme events in high
dimensional systems. The indicator is based on the aforementioned observations on the multiscale systems
but does not require a priori known fast-slow coordinates. More precisely we show that a small number of optimally time-dependent (OTD) modes \citep{otd}, obtained through a minimization principle and the history of the system state up to the current time instant,
allows for the description of the currently most unstable subspace  in
a dynamically
consistent fashion.
We show 
that the  linearized dynamics projected in this optimal, time-dependent subspace, can predict 
an upcoming rare extreme event. 
We note that simply computing the eigenvalues of the linearized dynamics is 
costly and, in many cases, the results are completely oblivious to transient instabilities (e.g. instabilities associated with non-normal dynamics, \citep{otd}).

In Section~\ref{sec:prelim}, we review the OTD modes and
introduce our indicator. We demonstrate the application of the indicator
on three examples: a low dimensional prototype system (Section~\ref{sec:bursting}),
a body forced Navier--Stokes equation (Section~\ref{sec:kolm})
and a modified nonlinear Schr\"odinger equation 
as a model for unidirectional water waves (Section~\ref{sec:mnls}).
The concluding remarks are presented in Section~\ref{sec:concl}.

\section{Preliminaries}\label{sec:prelim}
\subsection{Set-up}
Consider the general nonlinear system of ordinary differential equations (ODEs),
\begin{equation}
\dot{\vc u}=\vc F(\vc u),\quad \vc u\in\mathbb R^n,
\label{eq:ode}
\end{equation}
where the vector field $\vc F:\mathbb R^n\to \mathbb R^n$ is sufficiently smooth.
We denote the solutions of~\eqref{eq:ode} with the initial condition $\vc u_0$
at time $t_0$ by $\vc u(t;t_0,\vc u_0)=\varphi^t_{t_0}(\vc u_0)$ 
where $\varphi^t_{t_0}$ is the flow map.
Infinitesimal perturbations around an arbitrary trajectory $\vc u(t)$ satisfy the linear equation
\begin{equation}
\dot{\vc v}=\vc L_{\vc u}\vc v,\quad \vc v\in\mathbb R^n,
\label{eq:lin_ode}
\end{equation}
where $\vc L_{\vc u(t)}:=\bnabla\vc F(\vc u(t))$. For notational simplicity, we will write $\vc L$ instead of $\vc L_{\vc u}$.

For a given trajectory $\vc u(t;t_0,\vc u_0)$, there exists a two-parameter family of linear maps 
$\Phi_{t_0}^t(\vc u_0):\mathbb R^n\to \mathbb R^n$ such that the solutions of the linear equation~\eqref{eq:lin_ode} satisfy 
$\vc v(t;t_0,\vc v_0)=\Phi_{t_0}^t(\vc u_0)\vc v_0$ \citep{arnold_ode}. 
For notationally simplicity, we denote the solutions of the linear equation~\eqref{eq:lin_ode} by $\vc v(t)$
and write $\vc v(t)=\Phi_{t_0}^t\vc v_0$ along a given trajectory $\vc u(t)=\vc u(t;t_0,\vc u_0)$ of the nonlinear 
system~\eqref{eq:ode}.

In order to introduce the OTD modes, we will need the following definition.
\begin{defn}
A time-dependent $r$-dimensional subspace $E_r(t)$ of $\mathbb R^n$ 
is \emph{flow invariant} under the system~\eqref{eq:lin_ode} if, for a fixed 
initial time $t_0$, we have
\begin{equation}
\vc v(t)=\Phi_{t_0}^t\vc v_0\in E_r(t),\quad \forall \vc v_0\in E_r(t_0),\quad \forall t\geq t_0.
\end{equation}
\label{def:flowInv}
\end{defn}

\subsection{Optimally time-dependent modes}\label{sec:otd}
For $r=n$ in Definition~\ref{def:flowInv}, we have $E_n(t)=\mathbb R^n$ for all $t$ and therefore the space is trivially flow invariant.
Lower dimensional flow invariant subspaces can in principle be constructed as follows. Consider a
prescribed set of $r$ vectors $\{\vc v_1(t_0),\cdots,\vc v_r(t_0)\}$ spanning an $r$ dimensional subspace 
$E_{r}(t_0)$. For any later time $t>t_0$, let $\vc v_i(t)$ be the solutions of the liner equation~\eqref{eq:lin_ode} with the 
initial condition $\vc v_i(t_0)$ and define $E_r(t)=\mbox{span}\{\vc v_1(t),\cdots,\vc v_r(t)\}$. Since the map $\Phi_{t_0}^t$ is a bijection, the dimension of the linear subspace $E_r(t)$ is equal to $r$. Moreover, the subspaces
$E_r(t)$, constructed as such, are flow invariant by definition.
\begin{figure}
\centering
\includegraphics[width=.75\textwidth]{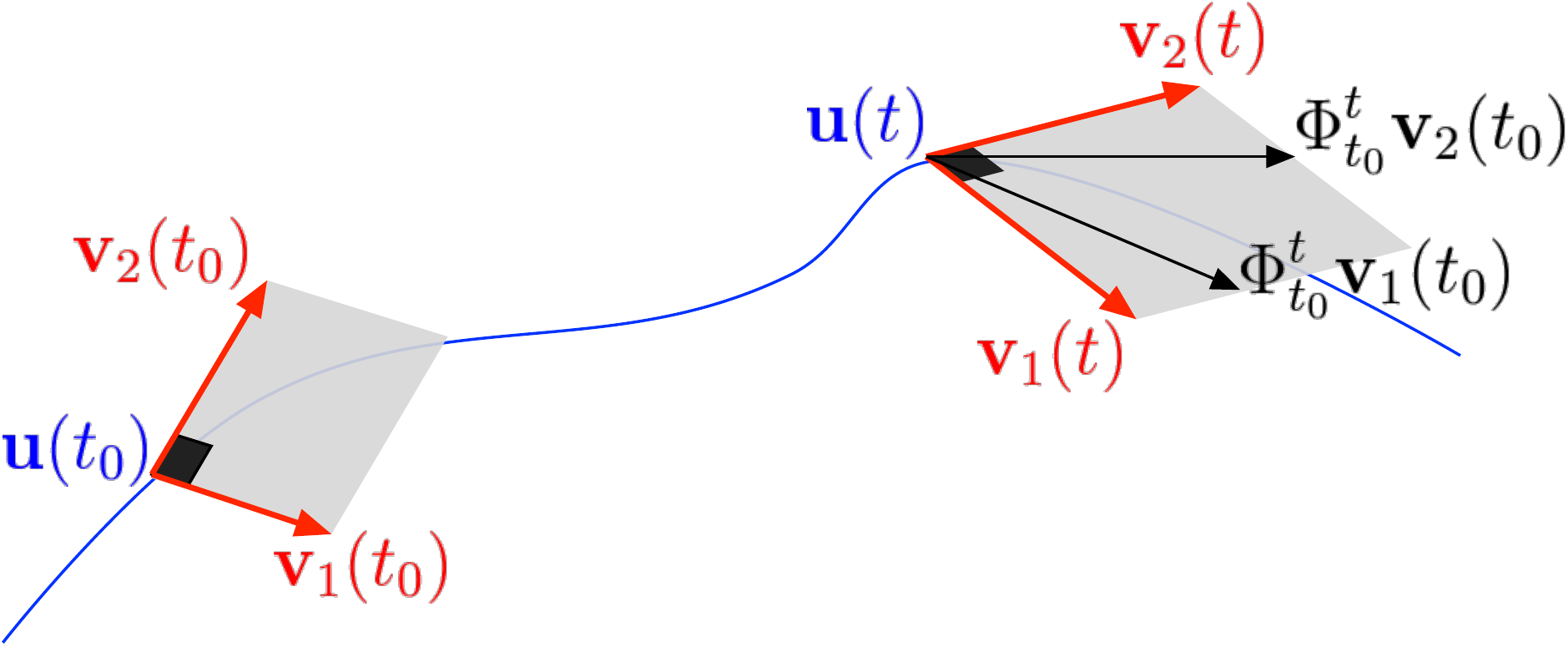}
\caption{An illustration of the OTD modes. The OTD modes $\vc v_i$ remain orthonormal 
for all times (the dark black 
squares mark right angles). While
differing from their images under the linear dynamics $\Phi_{t_0}^t$, the OTD modes span the same
subspace as their images.
}
\label{fig:schem_otd}
\end{figure}

This procedure is, however, known to be numerically unstable: typically the lengths of the vectors $\vc v_i$ grow 
exponentially fast and the angle between them vanishes rapidly. As a result, many numerical techniques have
been introduced to compute the flow invariant subspace in a numerically robust fashion 
(see, e.g.,~\citet{greene87,DE08}).

The OTD equations, introduced recently by~\citet{otd}, is a modification to the 
equation of variations~\eqref{eq:lin_ode} such that its solutions (the OTD modes)
remain orthonormal for all times, yet they span the same flow invariant subspaces as
the solutions of the equation of variations.

Here, we briefly review the OTD equations and the main properties of their solutions, referring
the interested reader to~\cite{otd} for details. The OTD equations read
\begin{equation}
\dot{\vc v}_i = \lins{\vc u} \vc v_i - \sum_{k=1}^{r}\langle \lins{\vc u}\vc v_i,\vc v_k\rangle \vc v_k,
\quad i\in\{1,2,\cdots, r\},
\label{eq:do_gen}
\end{equation}
where $\langle\cdot,\cdot\rangle$ denotes an appropriate inner product and $1\leq r\leq n$ is 
some prescribed integer.
Equations~\eqref{eq:do_gen}, together with the original system~\eqref{eq:ode},
form a set of $(r+1)$ coupled nonlinear differential equations for
vectors $\vc v_i\in\mathbb R^n$ and the state $\vc u$.
Note that without the summation term, the OTD equation~\eqref{eq:do_gen}
coincides with the equation of variations~\eqref{eq:lin_ode}.
The summation terms impose the constraint that the solutions $\vc v_i$
remain orthonormal with respect to the inner product $\langle\cdot,\cdot\rangle$.

We refer to the solutions $\vc v_i$ of the OTD equation as the \emph{OTD modes},
which have the following appealing properties.
\begin{enumerate}
\item The OTD equations preserve orthonormality: 
Let the initial condition
for the OTD equations~\eqref{eq:do_gen} be a set of orthonormal vectors
$\{\vc v_1(t_0),\vc v_2(t_0),\cdots,\vc v_r(t_0)\}$. Then
the solution $\{\vc v_1(t),\vc v_2(t),\cdots,\vc v_r(t)\}$ of the OTD equation
remains orthonormal for all times $t$~\citep[see][Theorem 2.1]{otd}.

\item The OTD modes span flow invariant subspaces: Define
\begin{equation}
E_r(t)={\rm span} \{\vc v_1(t),\vc v_2(t),\cdots,\vc v_r(t)\},
\end{equation}
with $\{\vc v_1(t),\vc v_2(t),\cdots,\vc v_r(t)\}$ being an orthonormal solution of 
the OTD equation~\eqref{eq:do_gen}. Then the subspaces $E_r(t)$ 
are flow invariant under the linear system~\eqref{eq:lin_ode}~\citep[see][Theorem 2.4]{otd}.

\item If $\vc u$ is a hyperbolic fixed point, the OTD modes generically converge to the
subspace spanned by the $r$ least stable eigenvectors of $\lins{\vc u}$ \citep[see][Theorem 2.3]{otd}.
\end{enumerate}
Figure~\ref{fig:schem_otd} illustrates the geometry of OTD modes for $r=2$.

\subsection{Reduction to the OTD modes}
Due to the flow invariance of the OTD modes, we can reduce the linear operator $\lins{\vc u}$
to the OTD subspaces $E_r(t)$ in a dynamically consistent fashion. 
More precisely, consider the solutions of the form $\vc v(t) = \vc V(t)\etab(t)$ where
$\vc V =\left[
\vc v_1 | \vc v_2 | \cdots | \vc v_r
\right]\in\mathbb R^{n\times r}$ 
is the time dependent matrix whose columns are the OTD modes
obtained from~\eqref{eq:do_gen}. 
The vector $\etab\in\mathbb R^r$ is the solution $\vc v$ expressed in the OTD basis.

Substituting $\vc v(t)=\vc V(t)\etab(t)$ in~\eqref{eq:lin_ode} yields the
\emph{reduced linear equation}
\begin{equation}
\dot{\etab}=\vc V^\dagger \vc L\vc V\etab. 
\label{eq:reduced_lin_ode}
\end{equation}
Conversely, if $\etab(t)$ solves the reduced equation~\eqref{eq:reduced_lin_ode}, then $\vc v(t)=\vc V(t)\etab(t)$ solves 
the full linear equation~\eqref{eq:lin_ode}~\citep[see][Theorem 2.4]{otd}.
We refer to the linear map $\vc L_r:\mathbb R^r\to\mathbb R^r$,
\begin{equation}
\vc L_r(t):=\vc V^\dagger(t) \vc L(t)\vc V(t),
\label{Lr}
\end{equation}
as the \emph{reduced linear operator}.

The reduced system~\eqref{eq:reduced_lin_ode} is a linear system of differential equations with
a time dependent stability matrix $\vc L_r(t)$. As a result, the eigenvalues of $\vc L_r$
may not be used to assess linear growth or decay of perturbations.
Instead we use the invariants of the symmetric part of 
$\vc L_r$ as an indicator. 

It follows from~\eqref{eq:reduced_lin_ode} that
\begin{equation}
\frac{1}{2}\frac{\id}{\id t}|\etab|^2=\langle\etab, \vc L_r\etab\rangle=\langle \etab, \vc S_r\etab\rangle,
\label{eq:insGrowth}
\end{equation}
where $\vc S_r$ denotes the symmetric part of the matrix $\vc L_r$, i.e.,
\begin{equation}
\vc S_r:=\frac{1}{2}\left[\vc L_r+\vc L_r^\dagger \right].
\label{eq:Sr}
\end{equation}
The eigenvalues $\lambda_1\geq\lambda_2\geq\cdots\geq\lambda_r$ 
of the symmetric tensor $\vc S_r\in\mathbb R^{r\times r}$, therefore,
measure the instantaneous linear growth (or decay) of perturbations
within the OTD subspace $E_r(t)$. Furthermore, 
the identity~\eqref{eq:insGrowth} implies the inequality
\begin{equation}
|\etab(t_0)|e^{\lambda_{min}(t-t_0)}\leq |\etab(t)|\leq |\etab(t_0)|e^{\lambda_{max}(t-t_0)},\quad \forall t\in[t_0,t_0+T],
\label{eq:finiteGrowth}
\end{equation}
for $T>0$ and $\lambda_{min}\leq \lambda_{max}$ defined as
\begin{equation}
\lambda_{min}:= \min_{\tau\in[t_0,t_0+T]} \lambda_r(\tau),\quad 
\lambda_{max}:=\max_{\tau\in[t_0,t_0+T]} \lambda_1(\tau).
\end{equation}
In particular, if $\lambda_{min}$ is positive, the perturbations within the OTD subspace $E_r(t_0)$ grow
exponentially fast over the time interval $[t_0,t]$.

Based on the above observation, we use the eigenvalue configuration
of the symmetric tensor $\vc S_r$ as the indicator for an upcoming burst.
In so doing, we assume that the OTD modes capture the most unstable flow invariant subspace
along a time-dependent trajectory. As pointed out in Section~\ref{sec:otd}, 
this has been proved by~\citet[][Theorem 2.3]{otd}
for hyperbolic fixed points, but remains an open problem for time-dependent trajectories.

In case the slow manifold is known as a graph over the slow variables, the 
connection between the largest eigenvalue $\lambda_1$ of the reduced symmetric 
tensor $\vc S_r$ and the instabilities transverse to the slow manifold can be made rigorous 
as shown by~\citet{hallersapsis10}. In practice, this graph is rarely known.

\section{Conceptual model}\label{sec:bursting}
For illustrative purposes, we construct a prototype system which has
simple dynamics with bursting episodes.
The system is described by the set of nonlinear ODEs,
\begin{align}
\dot x & = \alpha x + \omega y + \alpha x^2+2\omega xy + z^2\nonumber\\
\dot y & = -\omega x + \alpha y - \omega x^2+2\alpha xy\nonumber\\
\dot z & = -\lambda z -(\lambda+\beta)xz,
\label{eq:bursting}
\end{align}
where $\alpha,\omega,\lambda,\beta>0$ are constant parameters.
We define $\vc u=(x,y,z)$ and denote the right-hand-side of~\eqref{eq:bursting} by $\vc F(\vc u)$.
The plane $z=0$ is an invariant manifold containing the two fixed points 
\begin{equation*}
\vc u_1 = (0,0,0),\quad \vc u_2=(-1,0,0).
\end{equation*}
Linearizing around these fixed points, we obtain
\begin{equation}
\bnabla\vc F(\vc u_1)=
\begin{pmatrix}
\alpha & \omega & 0\\
-\omega & \alpha & 0\\
0 & 0 & -\lambda
\end{pmatrix},\quad
\bnabla\vc F(\vc u_2)=
\begin{pmatrix}
-\alpha & -\omega & 0\\
\omega & -\alpha & 0\\
0 & 0 & \beta
\end{pmatrix}.
\end{equation}
The plane $z=0$ is the linear unstable manifold $E^u$ of $\vc u_1$ 
corresponding to the eigenvalues $\alpha\pm\i\,\omega$.
The plane $z=0$ is also the linear stable manifold $E^s$ of the fixed point $\vc u_2$ 
with eigenvalues $-\alpha\pm\i\,\omega$. In the following, we set $\alpha=0.01$, $\omega=2\pi$ and
$\lambda=\beta =0.1$.
\begin{figure}
\centering
\subfigure[]{\includegraphics[width=.5\textwidth]{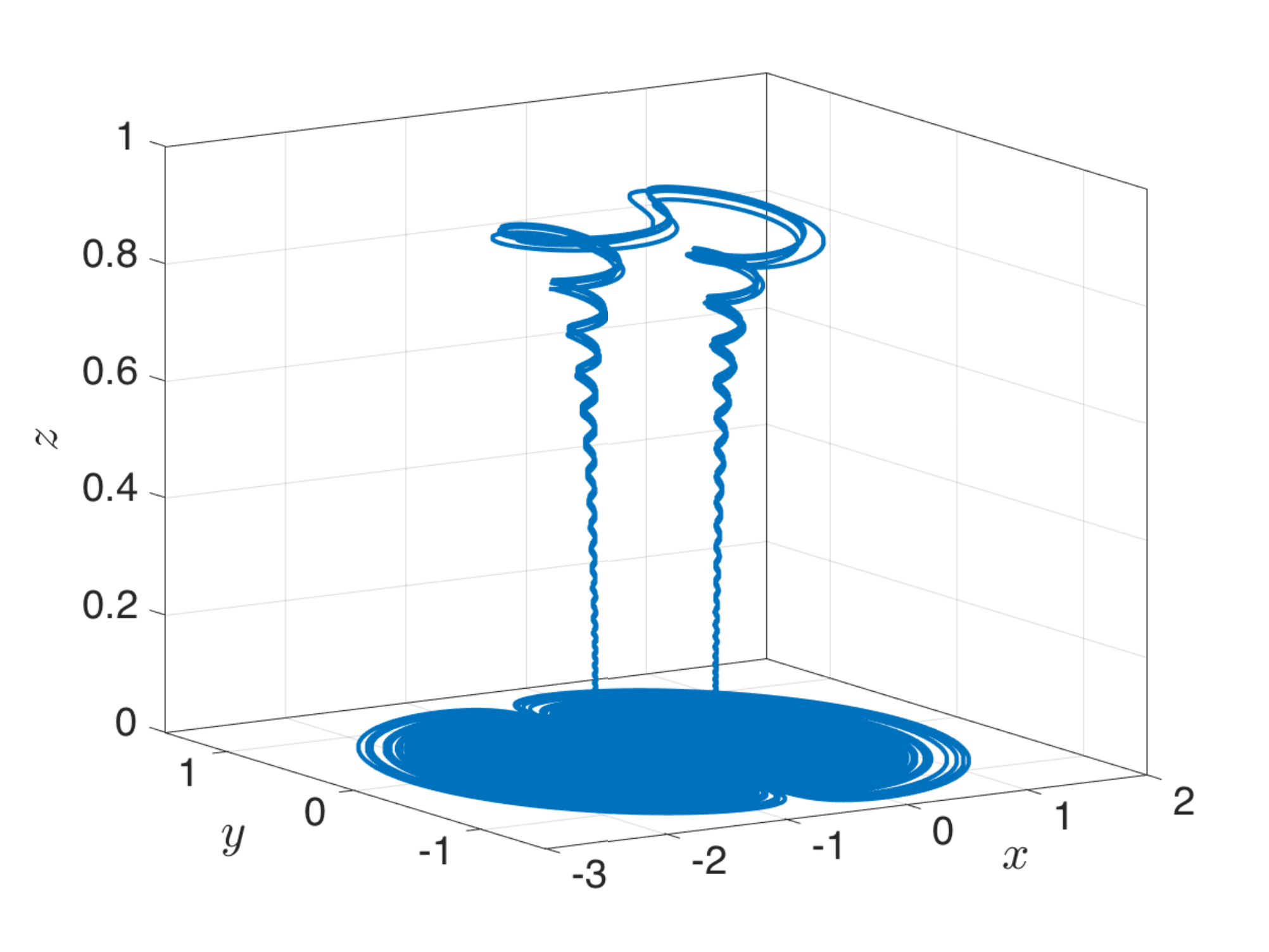}}
\subfigure[]{\includegraphics[width=.41\textwidth]{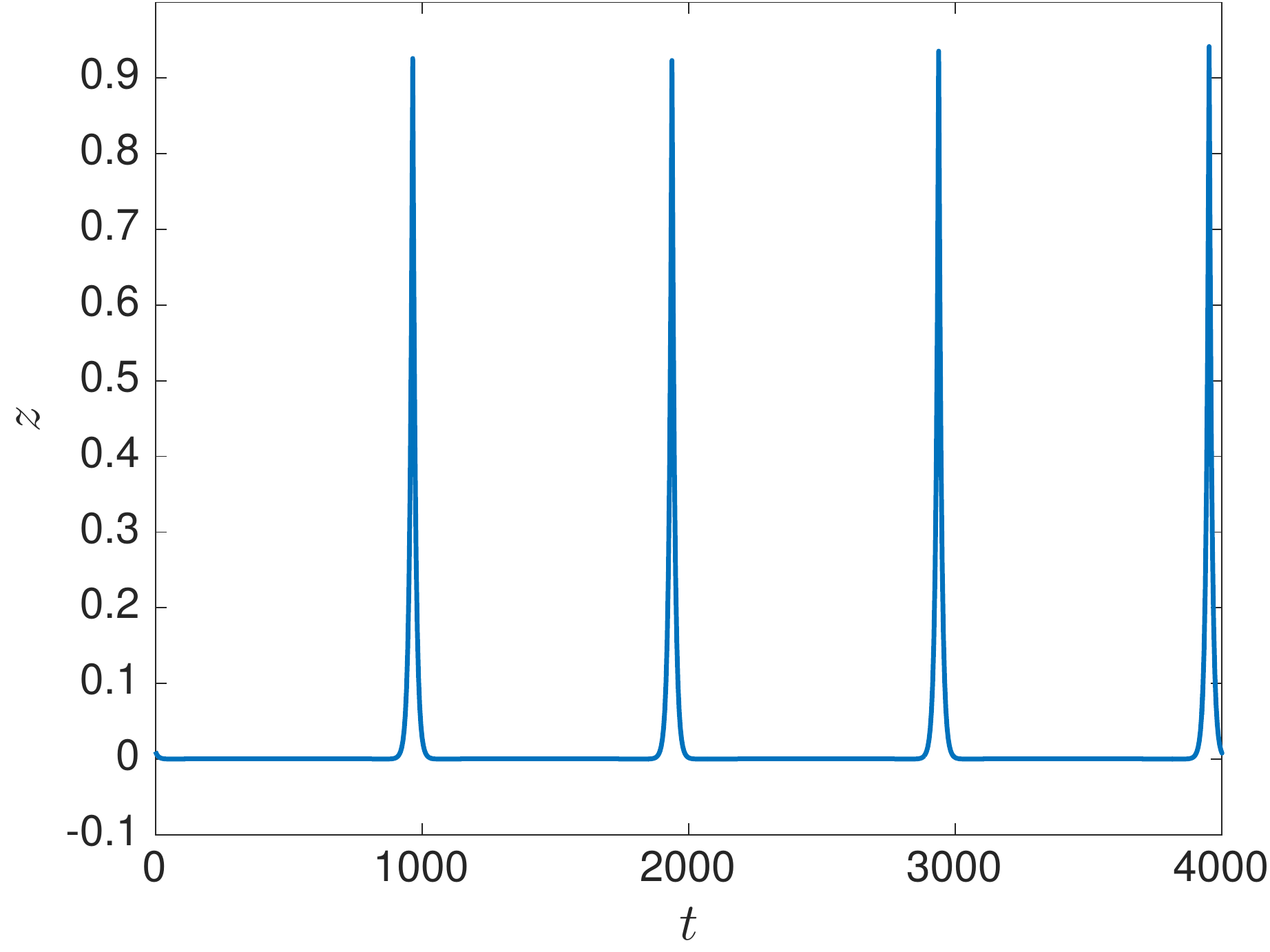}}
\caption{
A trajectory of the system~\eqref{eq:bursting} 
with parameters $\alpha=0.01$, $\omega=2\pi$, $\lambda=0.1$ and $\beta=0.1$.
The initial condition is $(0,0.01,0.01)$.
(a) The trajectory $\vc u(t)$ in the state space. 
(b) The time series of the $z$ component of the trajectory for $4\times 10^3$ time units.
}
\label{fig:bursting}
\end{figure}

Figure~\ref{fig:bursting} shows a trajectory of the system 
starting near the origin. Perturbations around the fixed point $\vc u_1$ spiral away from the origin due to the instability 
in the $z=0$ plane. Since $z=0$ is also the stable manifold of the fixed point $\vc u_2$, the 
perturbed trajectory is attracted towards $\vc u_2$. Due to the small stability exponent $\alpha=0.01$,  
this process takes place over a long period of time during which the $z$ component of the
trajectory stays small. Once close enough
to the fixed point $\vc u_2$, its unstable manifold repels the trajectory away from $z=0$ plane,
resulting in a rapid growth of the $z$ component.
Finally the trajectory is carried back to the fixed point $\vc u_1$ along the heteroclinic orbit 
connecting the two fixed points. The above process repeats once the trajectory is back
in the neighborhood of the origin $\vc u_1$.
\begin{figure}
\centering
\subfigure[]{\includegraphics[width=.42\textwidth]{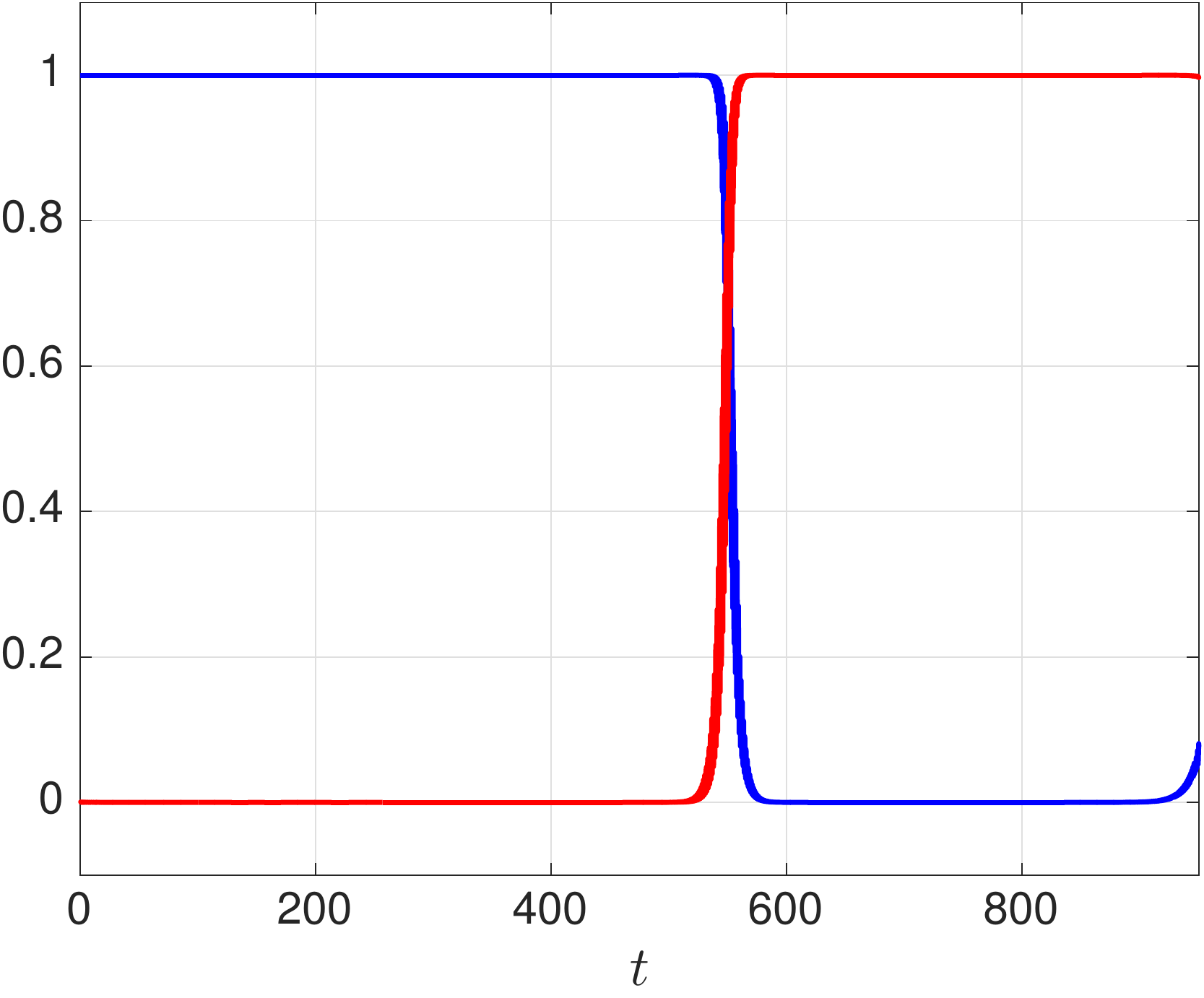}}
\subfigure[]{\includegraphics[width=.45\textwidth]{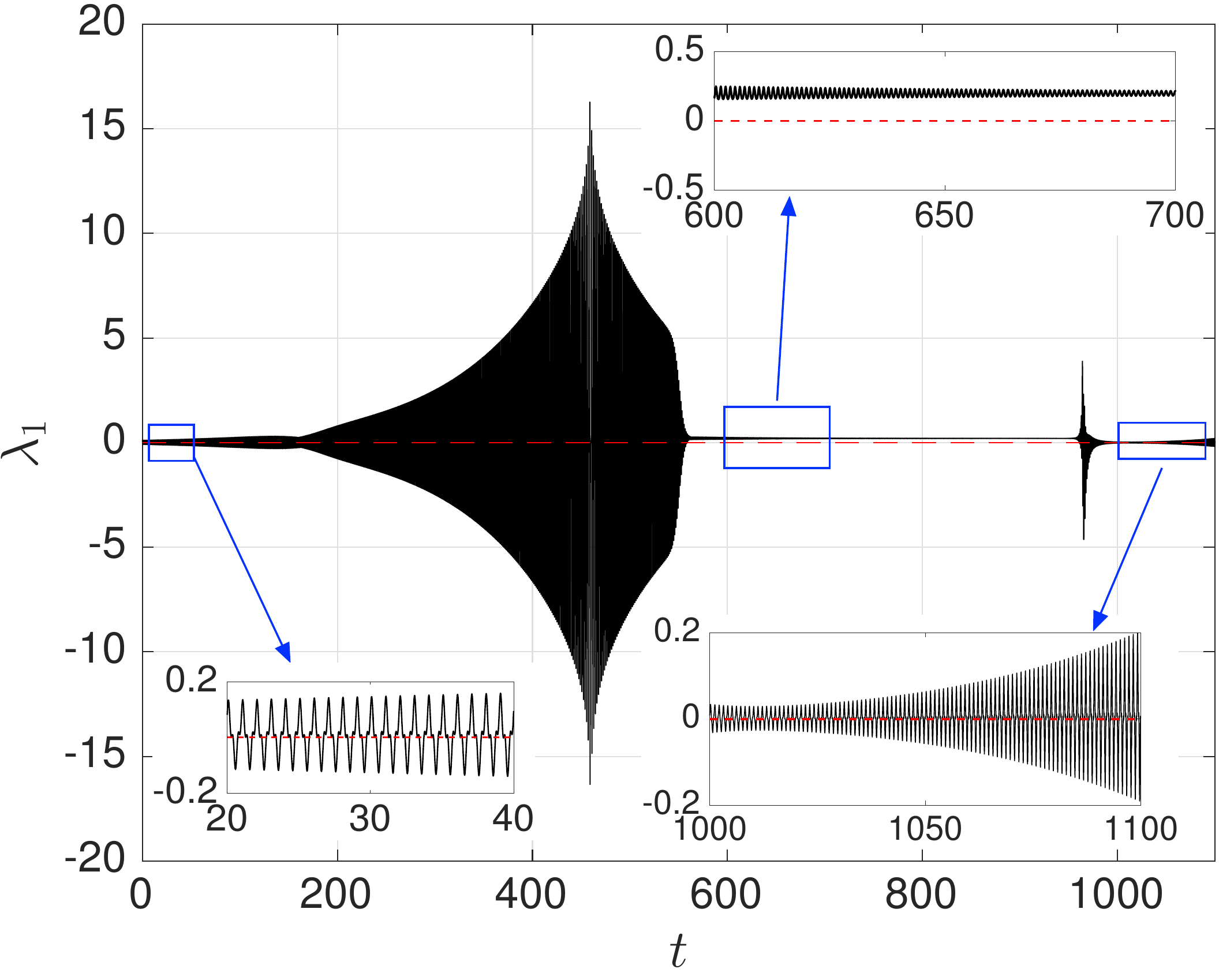}}
\caption{
(a) The evolution of $\sqrt{v_{1,1}^2+v_{1,2}^2}=\sqrt{1-v_{1,3}^2}$ (blue) and
$v_{1,3}$ (red) where $\vc v_1=(v_{1,1},v_{1,2},v_{1,3})$.
(b) The evolution of the eigenvalue $\lambda_1$ 
as a function of time. The dashed red line marks $\lambda_1=0$.
Three closeup views are shown in the insets.
}
\label{fig:bursting_otd}
\end{figure}

Now we investigate the ability of the OTD modes to capture the instability responsible for
the bursts. It is clear from the linearization that around the fixed point $\vc u_1$ the most unstable direction
is within the $x-y$ plane. Near fixed point $\vc u_2$ however the $z$-direction becomes the most 
unstable. We solve equation~\eqref{eq:bursting} together with 
the OTD equation~\eqref{eq:do_gen} with a single
OTD mode ($r=1$). We choose the initial conditions $\vc u=(0,0.01,0.01)^\top$
and $\vc v_1(0)=\frac{1}{\sqrt{2}}(1,1,0)^\top$.

Figure~\ref{fig:bursting_otd}(a) shows the evolution of $\sqrt{v_{1,1}^2+v_{1,2}^2}$ and $v_{1,3}$
where $v_{1,i}$ denote the components of $\vc v_1$, i.e., $\vc v_1=(v_{1,1},v_{1,2},v_{1,3})$.
For a long time, while the trajectory is spiraling away from $\vc u_1$, the $z$-component $v_{1,3}$
remains almost zero. As the trajectory moves towards the fixed point
$\vc u_2$, a sharp transition occurs around time $t=550$ where the OTD
mode $\vc v_1$ becomes almost orthogonal to the $x-y$ plane and aligns with the $z$ direction.
Note that this transition (at $t=550$) occurs well before the first burst (at $t=950$)
is observed (compare to figure~\ref{fig:bursting}(b)).

Figure~\ref{fig:bursting_otd}(b) shows the eigenvalue $\lambda_1$ of the reduced symmetric matrix 
$\vc S_r$ as a function of time. Since we only use one mode, the eigenvalue is trivial:
$\lambda_1(t)=\langle \vc v_1(t),\bnabla\vc F(\vc u(t))\vc v_1(t)\rangle$.
Over the initial $550$ time units, where the OTD mode $\vc v_1(t)$ is almost parallel to the
$x-y$ plane, the eigenvalue $\lambda_1$ oscillates rapidly around zero. As a result any 
instantaneous growth in the OTD subspace is rapidly counteracted by an instantaneous 
decay. After time $t=550$, when the OTD mode reorients orthogonally to the $x-y$ plane, 
the eigenvalue $\lambda_1$ becomes uniformly positive for a long period of time up until 
the bursting happens. This allows for persistent growth in the OTD subspace which aligns
with the $z$ axis in this period (cf. equation~\eqref{eq:finiteGrowth}). This instability persists up until the burst eventually happens 
around $t=960$. After the burst the eigenvalue $\lambda_1$ goes back to oscillating 
around zero.

\section{Turbulent fluid flow}\label{sec:kolm}
A ubiquitous feature of turbulent fluid flow is the intermittent bursts observed in the
time series of their measured quantities such as energy dissipation~\citep{frisch,faraz_adjoint}.
Even at moderate Reynolds numbers, the dimension of the turbulent attractors
are high. Best available estimates suggest that the attractor dimension scales almost linearly with the Reynold 
number~\cite{constantin88,robinson01,kevlahan07}. 
Moreover, no appropriate change of coordinates is available to 
decompose the system into slow and fast variables~\cite{berkooz93}. 
Consequently, intermittencies of turbulent fluid flow are particularly difficult to analyze 
and hence serve as a challenging example to test our indicator.

\subsection{Governing equations and preliminaries}
The two-dimensional Kolmogorov flow is the incompressible Navier--Stokes equations
\begin{equation}
\partial_t\vc u = -\vc u\cdot\bnabla\vc u - \bnabla p + \nu\Delta \vc u +\vc f, \quad \bnabla\cdot\vc u=0,
\label{eq:nse}
\end{equation}
with the sinusoidal forcing $\vc f=\sin(ny)\vc e_1$ where $\vc e_1=(1,0)^\top$
and $n$ is a positive integer~\cite{PlSiFi91}. 
The flow is defined on the torus $\vc x=(x,y)\in\mathbb T^2=[0,2\pi]\times [0,2\pi]$
(i.e., periodic boundary conditions). The solution is the time dependent pair of velocity field $\vc u(\vc x,t)$
and pressure $p(\vc x,t)$. The non-dimensional viscosity  $\nu$ is the inverse of
the Reynolds number, $\nu=1/Re$.

The energy $E$, energy dissipation $D$ and energy input $I$ of the system are defined as
\begin{align}
E(t)=\frac{1}{2L^2}\iint_{\mathbb T^2}|\vc u(\vc x,t)|^2\id \vc x,\quad & 
D(t)=\frac{\nu}{L^2}\iint_{\mathbb T^2}|\omega(\vc x,t)|^2,\nonumber\\
I(t)=\frac{1}{L^2}\iint_{\mathbb T^2}\vc u(\vc x,t)\cdot \vc f(\vc x,t)\id\vc x,
\label{eq:IDE}
\end{align}
where $L=2\pi$ is the size of the domain and $\omega$ is the vorticity field. 
One can show, from the Navier--Stokes equation~\eqref{eq:nse}, that these three quantities 
satisfy $\dot E=I-D$ along any trajectory.
\begin{figure}
\centering
\includegraphics[width=\textwidth]{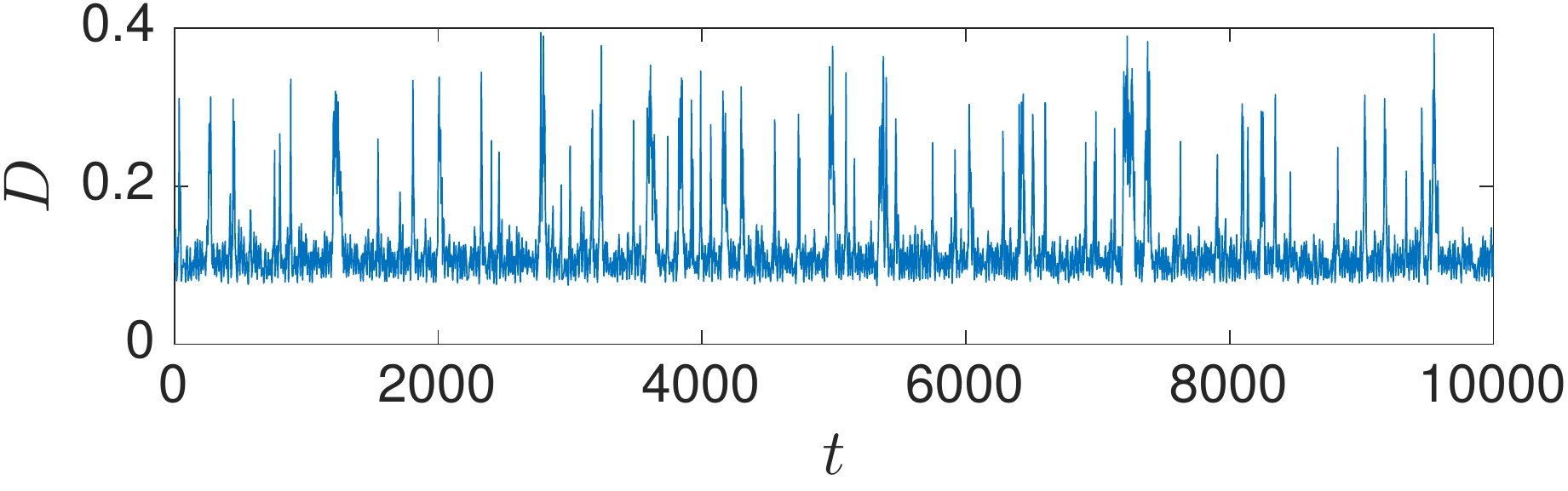}
\caption{Evolution of the energy dissipation $D$ along a trajectory of the Kolmogorov flow~\eqref{eq:nse}
with $n=4$ and $Re=40$.}
\label{fig:D_1001}
\end{figure}

The Kolmogorov flow has a \emph{laminar} solution,
\begin{equation}
\vc u_{lam}=\frac{Re}{n^2}\sin(ny)\vc e_1,
\label{eq:lam}
\end{equation}
which is asymptotically stable for forcing wave number $n=1$ and any Reynolds number $Re$
~\citep{marchioro86,foias2001}. For $n>1$ 
and sufficiently high $Re$, however, the laminar solution is unstable. 
Moreover, numerical evidence suggests that, for high enough Reynolds number
and $n>1$, the Kolmogorov flow is chaotic~\citep{PlSiFi91,CK13}.
Figure~\ref{fig:D_1001},
for instance, shows the evolution of the energy dissipation measured along a trajectory of the 
Kolmogorov flow with $n=4$ and $Re=40$. The energy dissipation mostly oscillates irregularly around
$D=0.1$ and never settles down to a regular pattern.
More interestingly, the energy dissipation exhibits intermittent, short-lived episodes of higher energy dissipation
that we wish to predict.

\subsection{OTD modes for the Kolmogorov flow}
In Section~\ref{sec:otd}, we introduced the OTD modes for ordinary differential equations. 
The OTD modes for partial differential equations (PDEs) are defined in a similar manner, although
more care should be exercised due to the infinite dimensionality of the system. 
In analogy with the ODEs, we define
\begin{equation}
\vc F(\vc u)=\mathbb P(-\vc u\cdot\bnabla\vc u + \nu\Delta \vc u +\vc f),
\end{equation}
where $\mathbb P$ denotes the projection onto space of divergence-free vector fields, $\bnabla\cdot\vc u=\vc 0$.
As opposed to the ODEs, where $\vc F$ is a vector field, here it is a nonlinear differential operator acting on 
functions $\vc u:\mathbb T^2\times\mathbb R \to \mathbb R^2$ that are sufficiently smooth.

We denote the linearization of $\vc F$ around the state $\vc u$ by $\lins{\vc u}$ whose action on
sufficiently smooth functions $\vc v:\mathbb T^2\times\mathbb R \to \mathbb R^2$ is given by
\begin{equation}
\lins{\vc u}\vc v:=\mathbb P(-\vc u\cdot \bnabla\vc v-\vc v\cdot \bnabla\vc u +\nu \Delta\vc v).
\label{eq:L_kolm}
\end{equation}

The OTD modes $\{\vc v_1,\vc v_2,\cdots,\vc v_r \}$ then satisfy the set of coupled PDEs
\begin{equation}
\frac{\partial \vc v_i}{\partial t}=\lins{\vc u}\vc v_i-
\sum_{j=1}^{r}\langle\lins{\vc u}\vc v_i,\vc v_j\rangle \vc v_j,\quad i\in\{1,2,\cdots,r\},
\label{eq:otd_pde}
\end{equation}
where $\langle\cdot,\cdot\rangle$ denotes some appropriate inner product. Here, we use the $L^2$ inner product
\begin{equation}
\langle\vc v,\vc w\rangle := \iint_{\mathbb T^2}\vc v(\vc x,t)\cdot \vc w(\vc x,t)\id \vc x.
\end{equation}
We integrate equations~\eqref{eq:otd_pde} with initial conditions
\begin{equation}
\vc v_k(\vc x,0)=\frac{1}{\pi\sqrt 2}
\begin{pmatrix}
\sin(ky)\\
0
\end{pmatrix},\quad k=1,2,\cdots,r\, ,
\label{eq:ic_otd}
\end{equation}
which are divergence free, mutually orthogonal and have unit $L^2$ norm.

The restriction of the infinite-dimensional operator $\lins{\vc u}$ to the time-dependent OTD 
subspace $\{\vc v_k \}_{1\leq k\leq r}$ is a reduced finite-dimensional linear operator $\vc L_r$. 
In the OTD basis, the reduced operator $\vc L_r$ is given by the $r\times r$ matrix
whose entries are given by
\begin{equation}
[\vc L_r]_{ij}=\langle \vc v_i,\lins{\vc u}{\vc v_j}\rangle,\quad i,j\in\{1,2,\cdots, r\}.
\end{equation}

Although the linear operator~\eqref{eq:L_kolm} acts on an infinite dimensional function space, 
the reduced operator $\vc L_r$ is a finite dimensional linear map whose 
symmetric part $\vc S_r$ is defined by~\eqref{eq:Sr}.

We numerically integrate the Kolmogorov equation~\eqref{eq:nse} and its associated OTD equations~\eqref{eq:otd_pde}.
To evaluate the right hand sides of the equations, we use a standard pseudo-spectral scheme with 2/3 
dealiasing~\cite{peyret13}. Unless stated otherwise, $128\times 128$ Fourier modes are used for the 
simulations reported below.
For the time integration, we use the Runge--Kutta scheme \texttt{RK5(4)}
of~\citet{RK45} with relative and absolute error tolerances set to $10^{-5}$.

\subsection{Asymptotically stable regime}
As mentioned earlier, for the forcing wavenumber $n=1$, the laminar solution~\eqref{eq:lam} of the Kolmogorov 
equation~\eqref{eq:nse} is asymptotically stable at any $Re$ number. Moreover, the laminar solution is also the 
global attractor~\citep{foias2001}. This regime is not our primary interest. It, however, does help illustrate the evolution 
of the OTD modes in an unambiguous setting. 
\begin{figure}
\centering
\includegraphics[width=.8\textwidth]{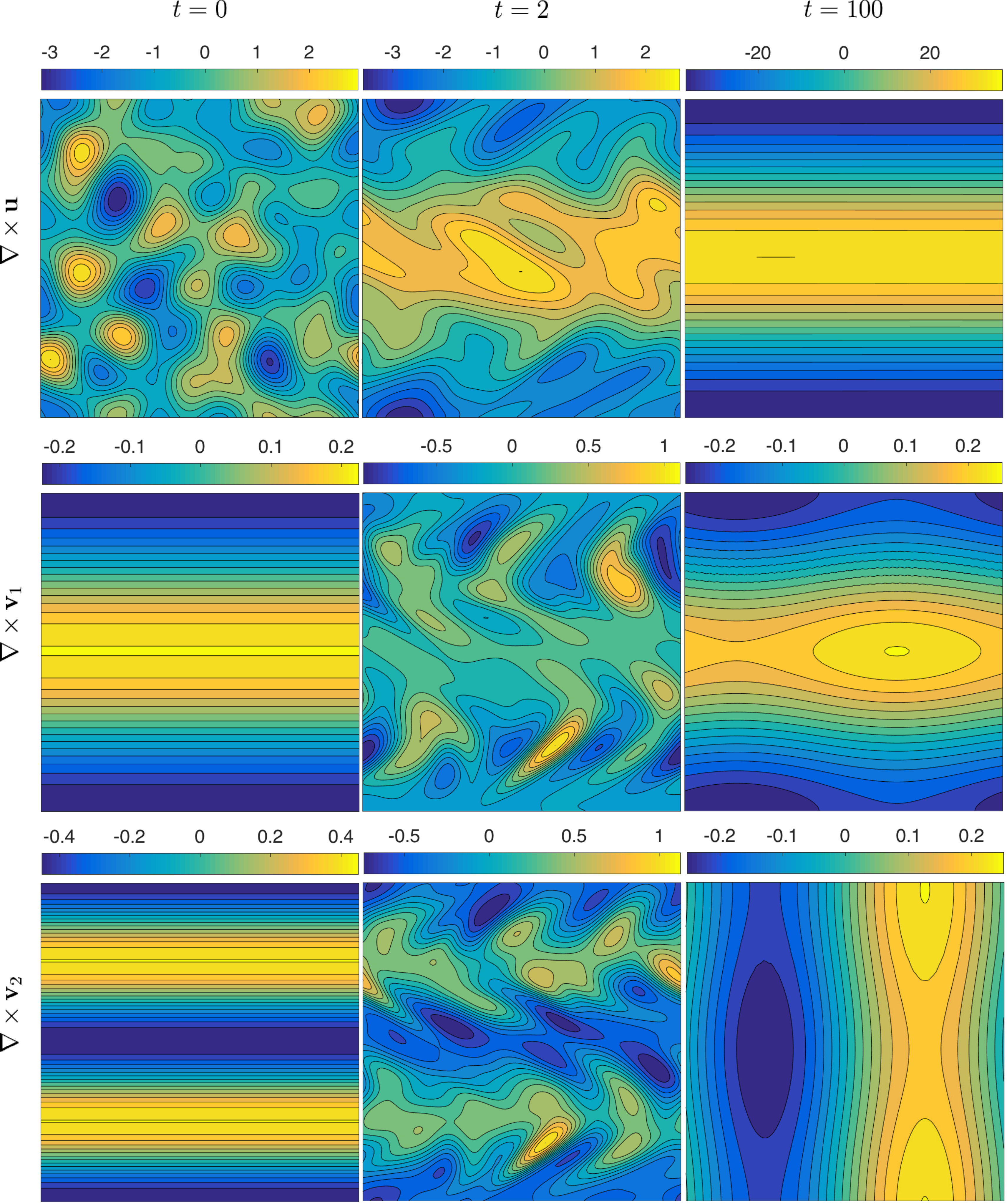}
\caption{The Kolmogorov flow in the asymptotically stable regime with $n=1$ and $Re=40$. 
Top row: The vorticity field at $t=0$, $2$ and $100$.
Middle row: curl of the first OTD mode $\vc v_1$ at $t=0$, $2$ and $100$.
Bottom row: curl of the second OTD mode $\vc v_2$ at $t=0$, $2$ and $100$.
The colors correspond to the only non-zero component of the curls.
All panels show the entire domain $[0,2\pi]\times [0,2\pi]$.
}
\label{fig:n1}
\end{figure}

We numerically solve the Kolmogorov equation and its associated OTD equations with $r=2$.
The state $\vc u$ is initially random in phase with an exponentially decaying energy spectrum.
The initial conditions for the OTD modes are given in~\eqref{eq:ic_otd}.
Figure~\ref{fig:n1}
shows the initial condition and the evolution of the state $\vc u$ and the OTD modes
$\vc v_1$ and $\vc v_2$ at select time instances.

\begin{figure}
\centering
\includegraphics[width=.45\textwidth]{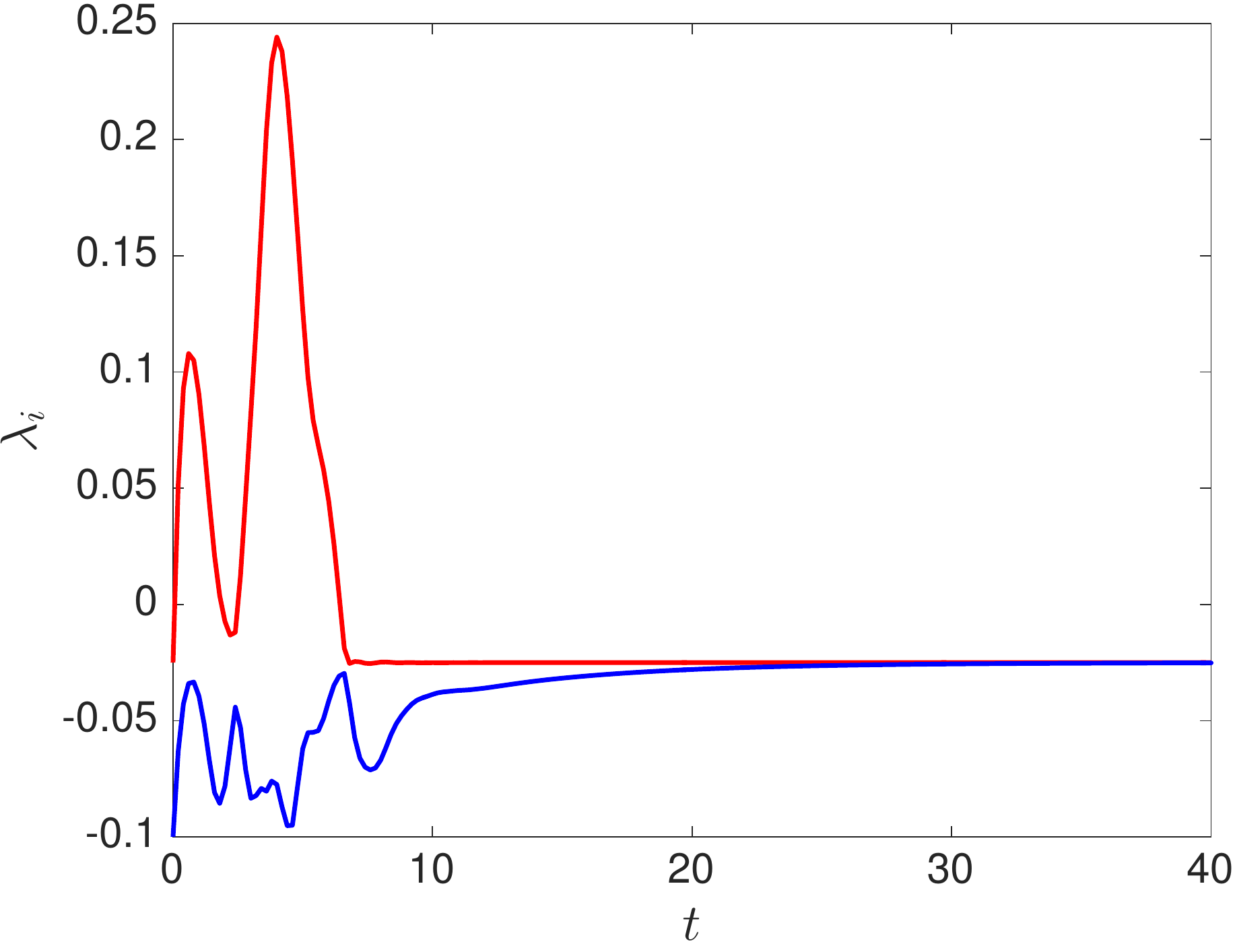}
\caption{Eigenvalues of the symmetric matrix $\vc S_2$ along a trajectory 
of the Kolmogorov flow in the asymptotically stable regime: $n=1$ and $Re=40$.}
\label{fig:n1_eig}
\end{figure}

The eigenvalues of the symmetric tensor $\vc S_2$ are shown in figure~\ref{fig:n1_eig}. 
As the flow evolves towards the laminar solution, the eigenvalues of $\vc S_2$ 
oscillate before they converge to their asymptotic value of $-0.025$. One of the eigenvalues
assumes positive values during this transition, signaling perturbations that can instantaneously grow.
Since the laminar solution is the global attractor, the instantaneous growth cannot be
sustained and decays eventually.
As the state $\vc u(t)$ converges to the laminar solution~\eqref{eq:lam}, the 
OTD modes $\vc v_1$ and $\vc v_2$ converge to the least stable eigenspace
of the linear operator~\eqref{eq:L_kolm} corresponding to
eigenvalue $-0.025$ whose algebraic and 
geometric multiplicity happens to be equal to $2$.

\subsection{Chaotic regime}
We turn now to a set of parameters for which the Kolmogorov flow is chaotic. 
Numerical evidence suggests that, for $n=4$ and $Re=40$, the Kolmogorov flow has
a strange attractor~\citep{CK13}. More importantly, the energy dissipation $D$
exhibits an intermittent behavior along the trajectories on the strange attractor (see figure~\ref{fig:D_1001}).

\begin{figure}
\centering
\subfigure[]{\includegraphics[width=.4\textwidth]{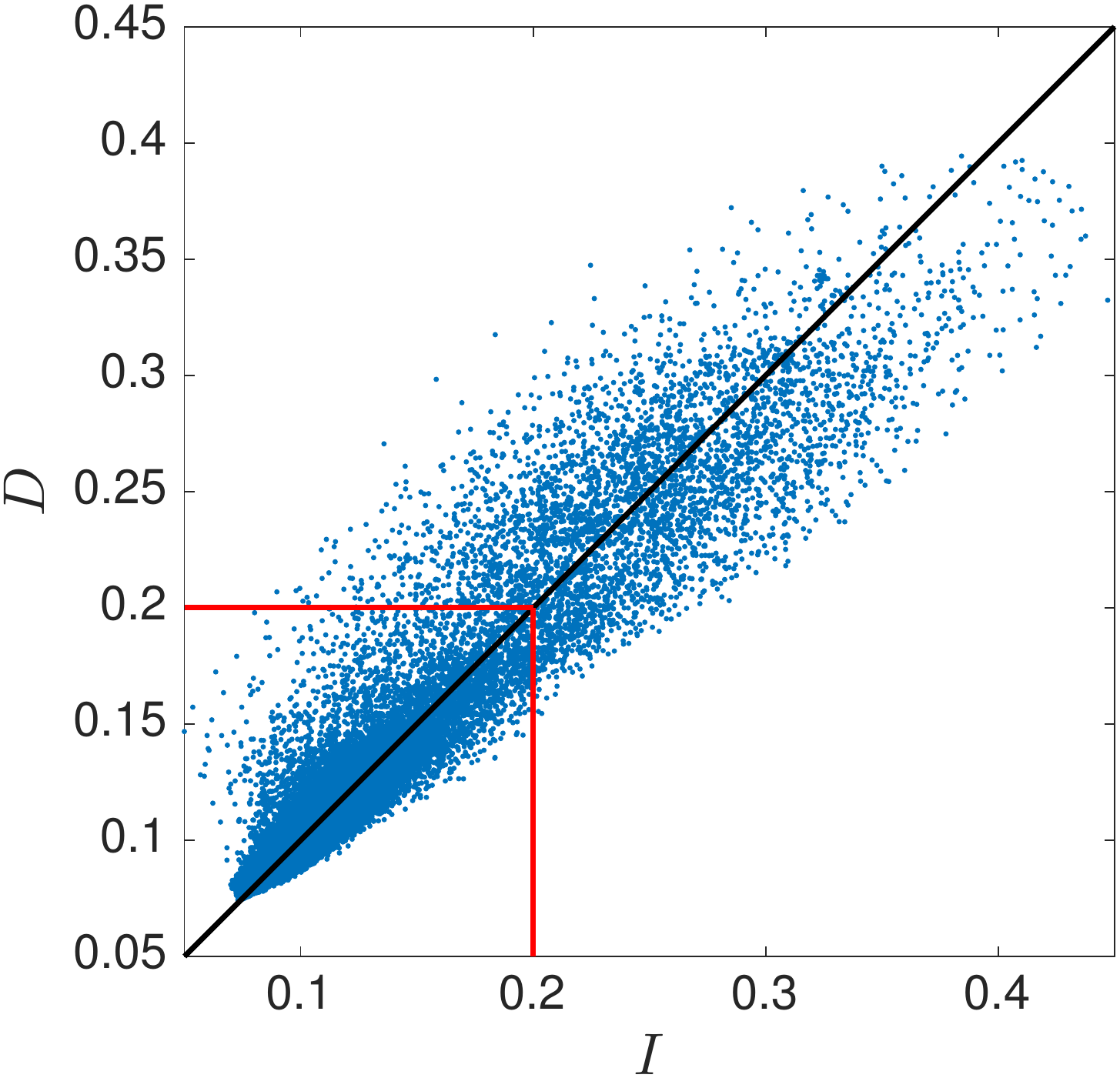}}
\subfigure[]{\includegraphics[width=.5\textwidth]{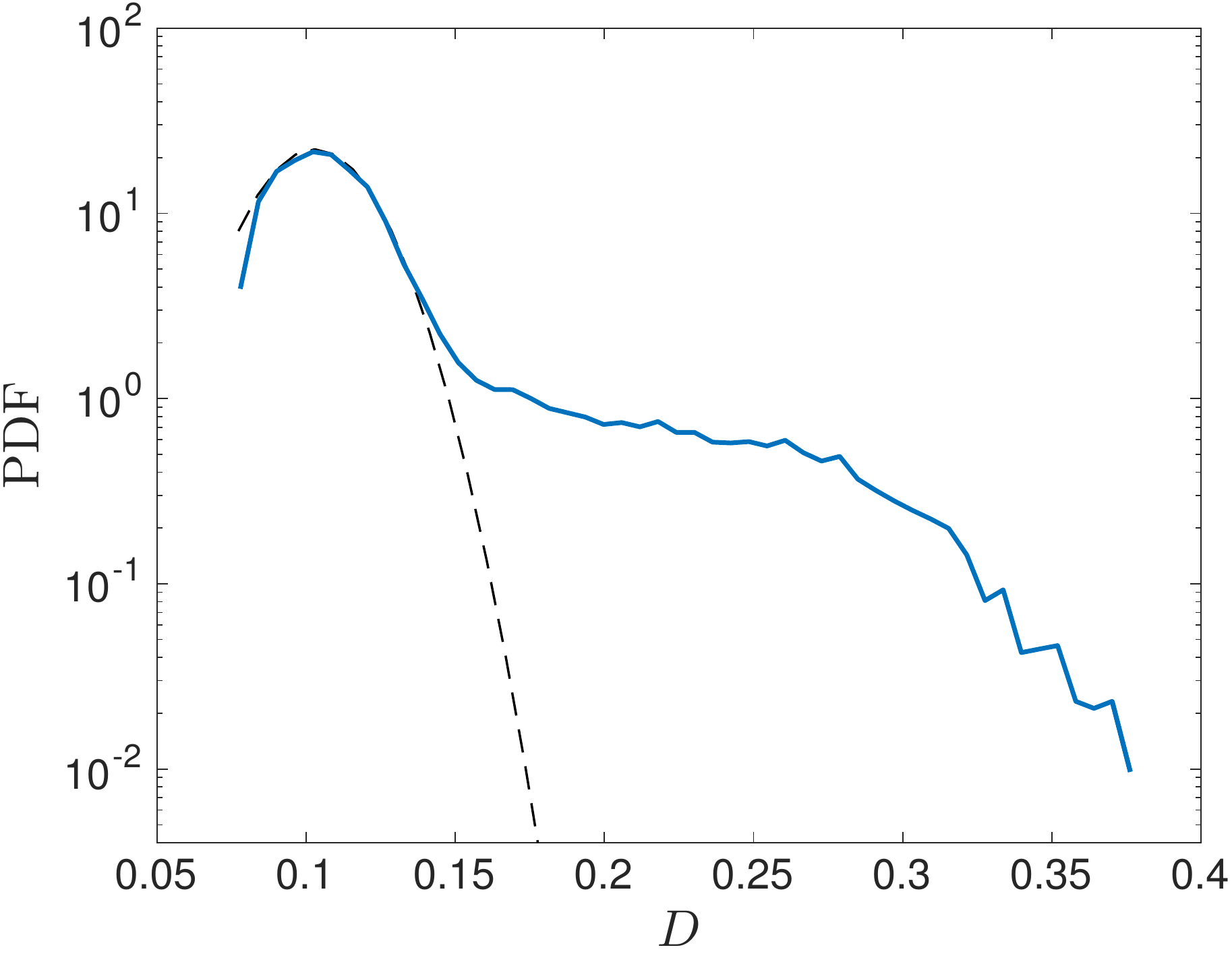}}
\caption{
(a) Energy input $I$ versus energy dissipation $D$ shown for a long turbulent trajectory. The dots correspond to
$5\times 10^4$ time instances each $0.2$ time units apart. The trajectory spends approximately $91.8\%$ of its lifetime 
inside the red box. The black line is the diagonal $I=D$.
(b) The probability density function (PDF) of the energy dissipation. The dashed black line marks the
PDF of a Gaussian districution with mean $0.103$ and standard deviation $0.018$.
}
\label{fig:300_I}
\end{figure}

Figure~\ref{fig:300_I} shows the energy input $I$ versus the energy dissipation $D$ for a long turbulent trajectory.
During the evolution, the energy input and dissipation assume smaller values and are very close to each other
sitting near the diagonal. The Kolmogorov flow is driven by the external forcing $\vc f$ such that 
growth in the energy input $I$
corresponds to the alignment of the velocity field $\vc u$ and the forcing (see equation~\eqref{eq:IDE}).
This alignment leads to an abrupt increase in the energy input $I$. Consequently, the energy dissipation
also increases bringing the trajectory back to the statistically stationary background.

Based on this observation, one may argue that the growth of the perturbations aligning with
the forcing should signal an upcoming burst in the energy input (and consequently the energy dissipation).
The instantaneous growth of such a perturbation is measured by $\langle \vc f,\lins{\vc u}\vc f\rangle$
(cf. equation~\eqref{eq:insGrowth}). For any divergence free velocity field $\vc u(t)$ with zero mean, however, 
a straightforward calculation yields $\langle \vc f,\lins{\vc u}\vc f\rangle\simeq -7.896$.
This seemingly paradoxical result is the consequence of the fact that the forcing $\vc f$
is not a flow invariant subspace and, as such, the instantaneously negative
value of $\langle \vc f,\lins{\vc u}\vc f\rangle$ does not imply decay over finite
time intervals. The OTD subspaces, in contrast, are flow invariant and therefore
a projection onto them is dynamically meaningful.

\begin{figure}
\centering
\includegraphics[width=.8\textwidth]{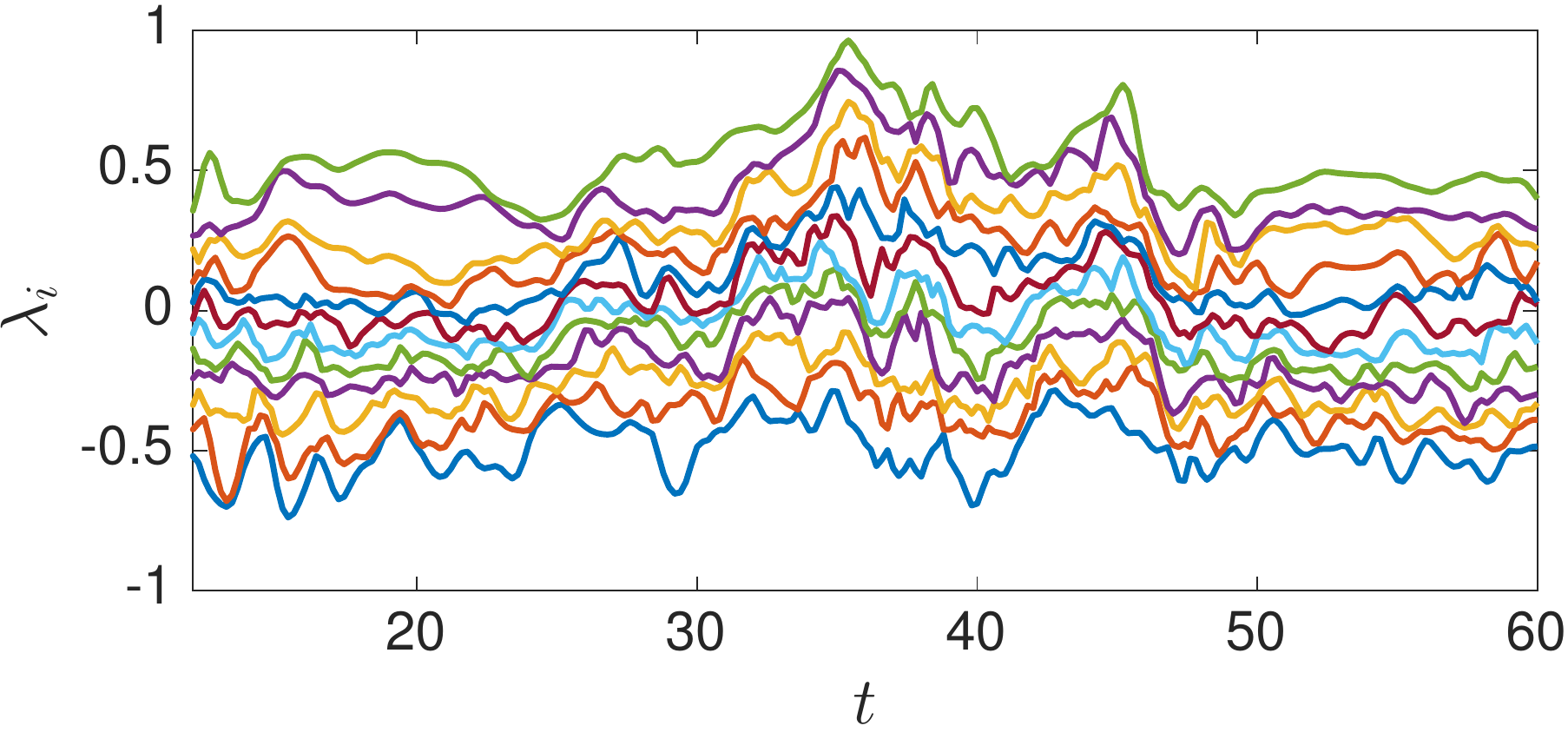}
\caption{Evolution of the eigenvalues $\lambda_1\geq \lambda_2\geq \cdots\geq \lambda_{12}$
of the reduced symmetric tensor $\vc S_{12}$ along a chaotic trajectory of the Kolmogorov flow.}
\label{fig:R40_eigs_L12s}
\end{figure}

The evolution of the eigenvalues of the symmetric tensor $\vc S_{12}$ along a turbulent trajectory
are shown in figure~\ref{fig:R40_eigs_L12s}. The first four eigenvalues are 
positive for all $t$ in this time window, signaling the very unstable nature of the flow.

\begin{figure}
\centering
\includegraphics[width=.8\textwidth]{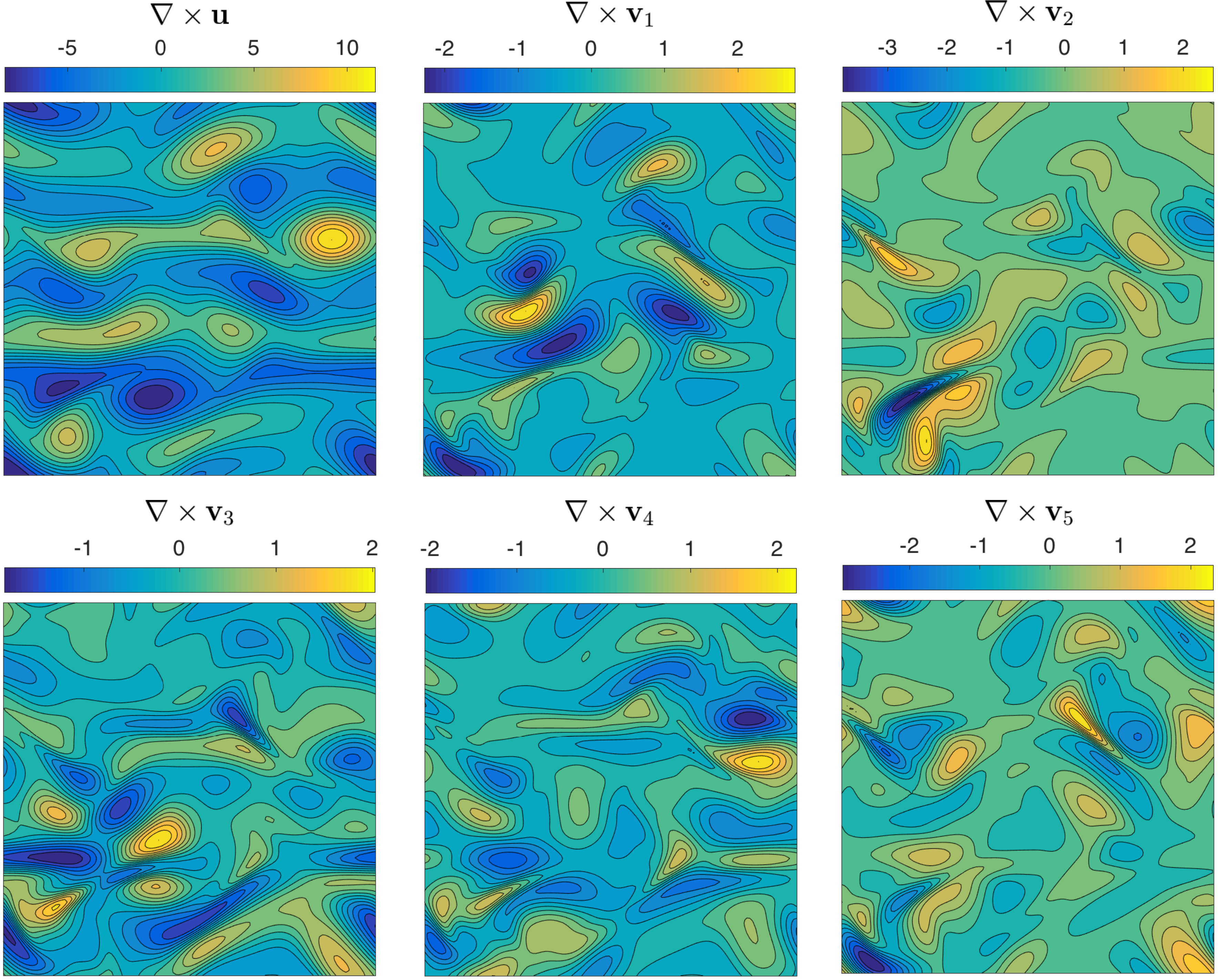}
\caption{Snapshots of the curl of 
$\vc u$, $\vc v_1$, $\vc v_2$, $\vc v_{3}$, $\vc v_{4}$ and $\vc v_{5}$
(from top left to bottom right, respectively) at time $t=34.6$.
The colors correspond to the only non-zero component of the curls.
All panels show the entire domain $[0,2\pi]\times [0,2\pi]$.}
\label{fig:otd_modes}
\end{figure}

\begin{figure}
\centering
\includegraphics[width=.48\textwidth]{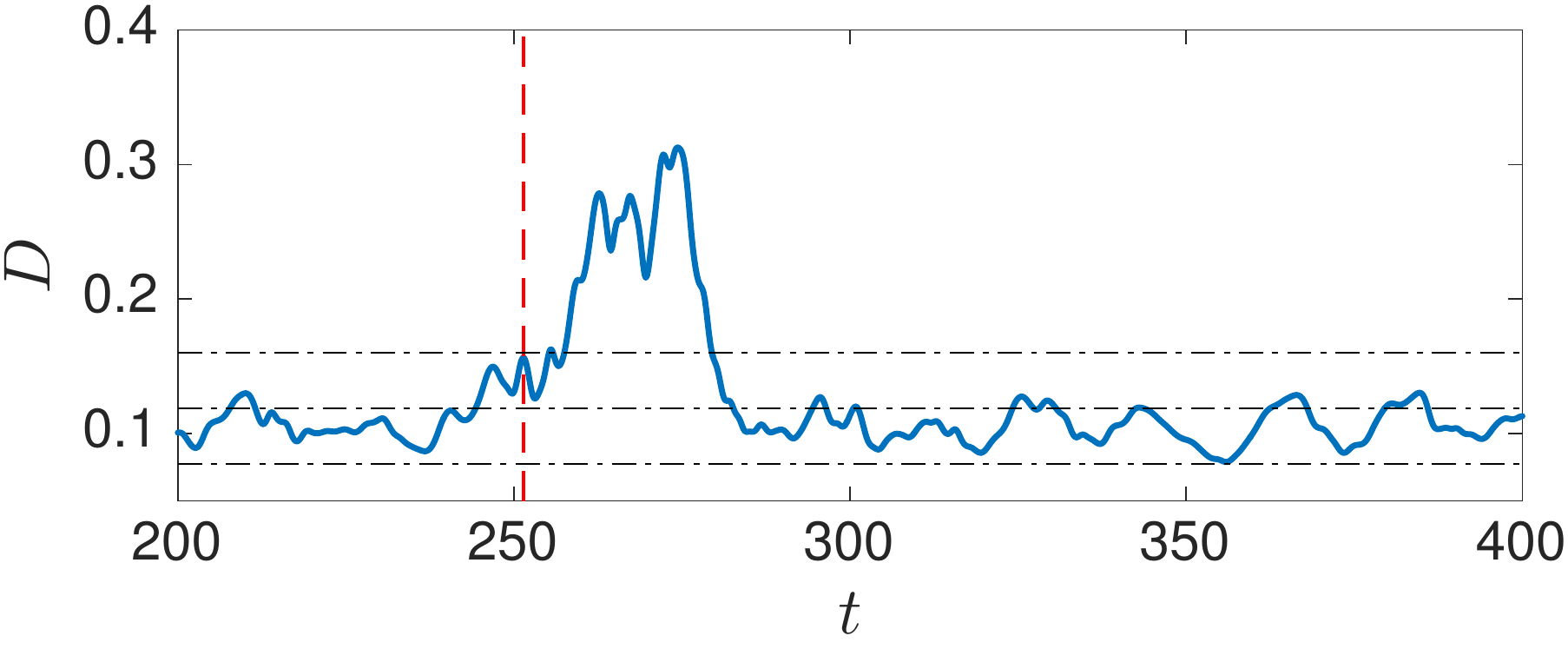}
\includegraphics[width=.48\textwidth]{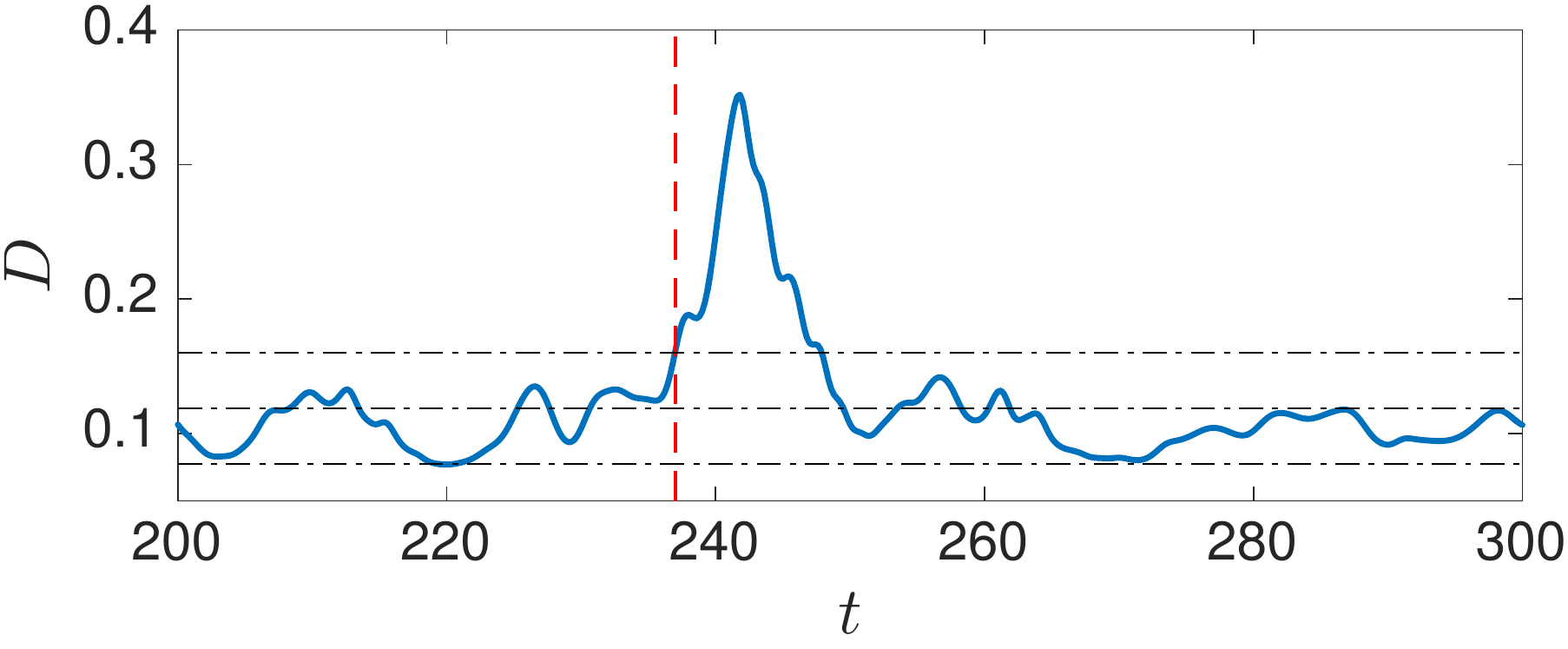}\\
\includegraphics[width=.48\textwidth]{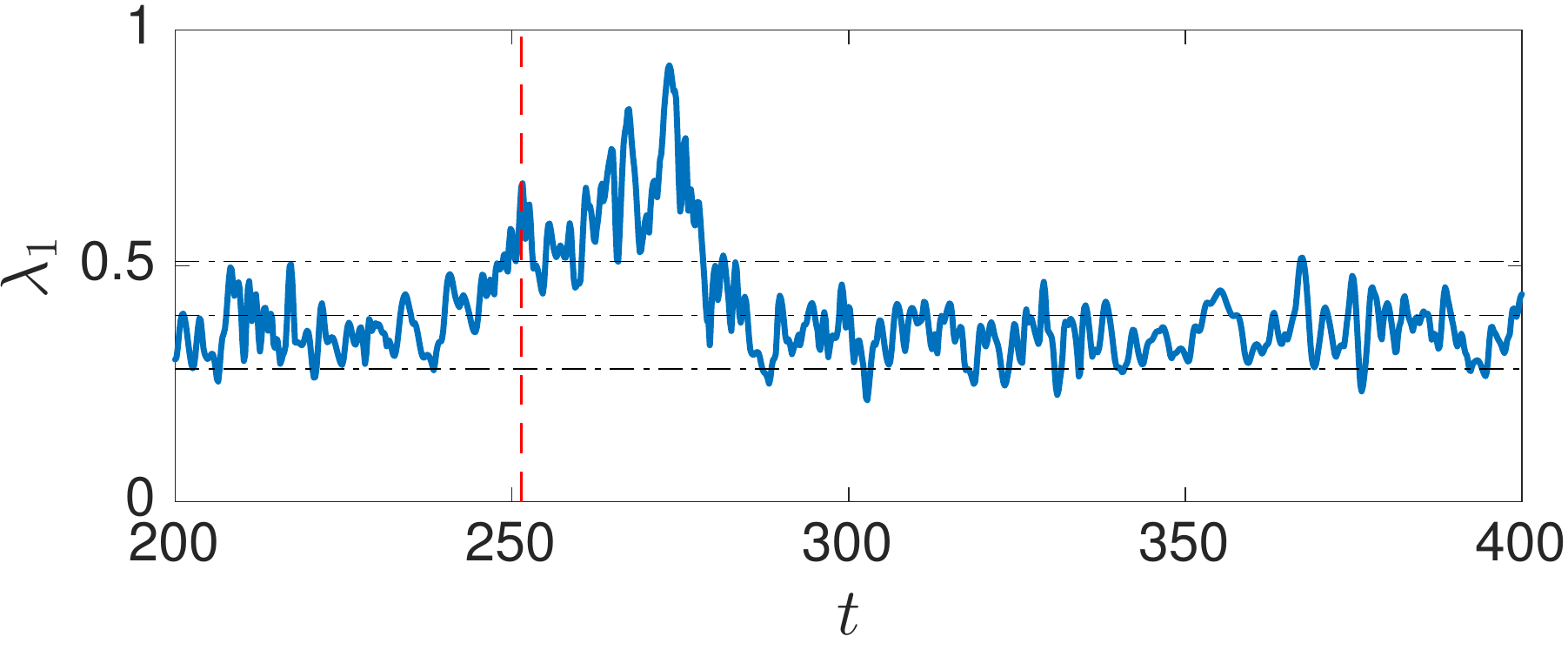}
\includegraphics[width=.48\textwidth]{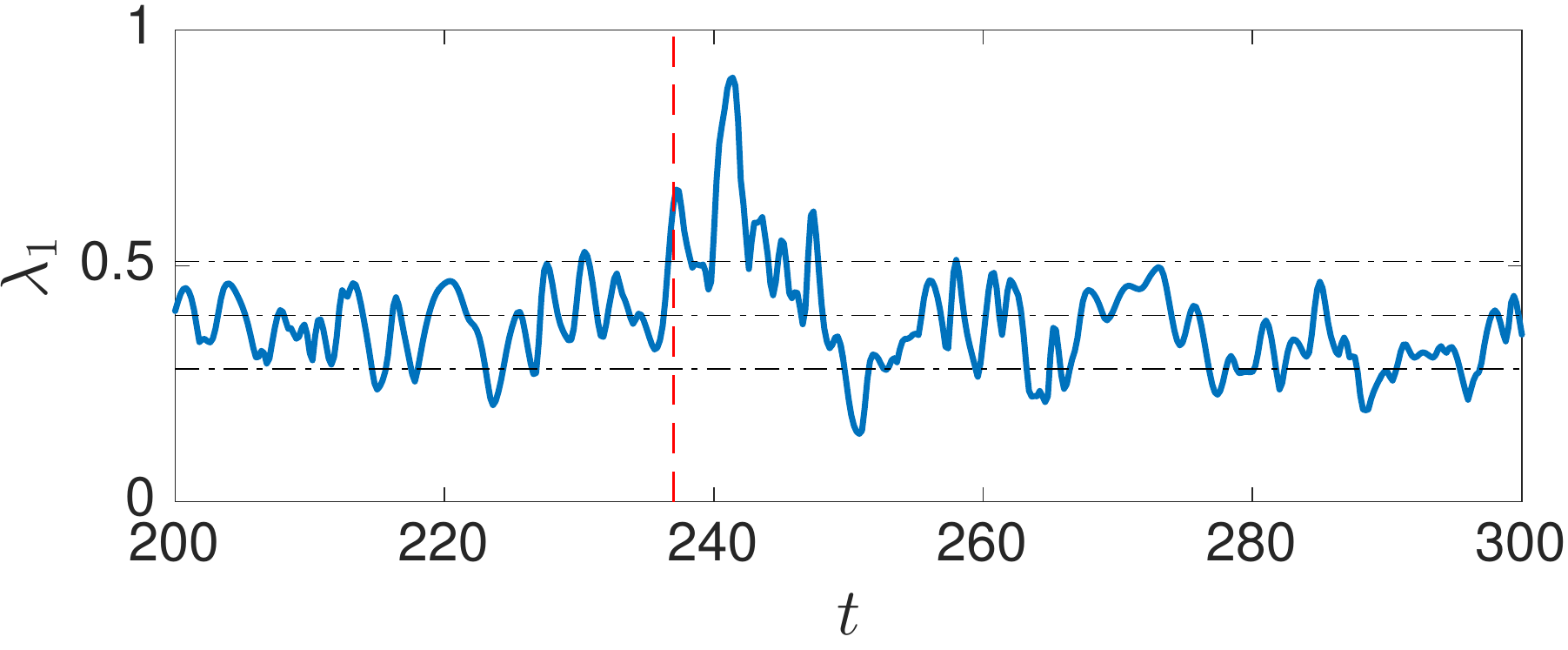}
\caption{Evolution of the energy dissipation $D$ and the eigenvalue $\lambda_1$
along two different trajectories (each column corresponds to a separate trajectory).
The horizontal dashed lines mark the mean, 
the mean plus one standard deviation and the mean minus one standard deviation
of the corresponding quantity.
}
\label{fig:1001_D_l1}
\end{figure}

Figure~\ref{fig:otd_modes} shows select OTD modes at time $t=34.6$ right before 
a burst in the energy dissipation occurs. The modes themselves do not exhibit a distinguished
structural feature that could suggest an immediate connection to the burst. 
We notice, however, that the largest eigenvalue $\lambda_1$ 
of the symmetric part of the reduced linear operator $\vc L_r$ 
increases significantly just before the bursting (see figure~\ref{fig:1001_D_l1}) 
while the energy dissipation is within one standard deviation from its mean value at that time.
Since the eigenvalue tends to oscillate rapidly and irregularly, mere visual 
inspection does not yield a satisfactory conclusion. 
In the next section, we quantify the correlation between the eigenvalue $\lambda_1$
and the energy dissipation $D$ through conditional statistics. 

\subsection{Conditional statistics}\label{sec:pdf}
In order to quantify the predictive power of the eigenvalues of reduced symmetric matrix
$\vc S_r$, we use Bayesian statistics~\citep{gelman14}. 
First, for a given scalar quantity $q(t)$, we define
\begin{equation}
\bar q(t;t_i,t_f)=\max_{\tau \in[t+t_i,t+t_f]}q(\tau),
\label{eq:Jq}
\end{equation}
where $0<t_i<t_f$ are prescribed numbers. At any time $t$, the quantity $\bar q(t;t_i,t_f)$
equals the maximum value of $q$ over a future time interval $[t+t_i,t+t_f]$.
For notational simplicity, we use the shorthand notation $\bar q(t)$ for $\bar q(t;t_i,t_f)$.

We would like to quantify the predictive power of a given indicator
$\alpha(t)$. In particular, we would like to assess whether large peaks of the indicator
$\alpha(t)$ at a time $t$ coincide with large values of the observable $q$ over 
a future time interval $[t+t_i,t+t_f]$.
\begin{figure}
\subfigure[]{\includegraphics[width=.3\textwidth]{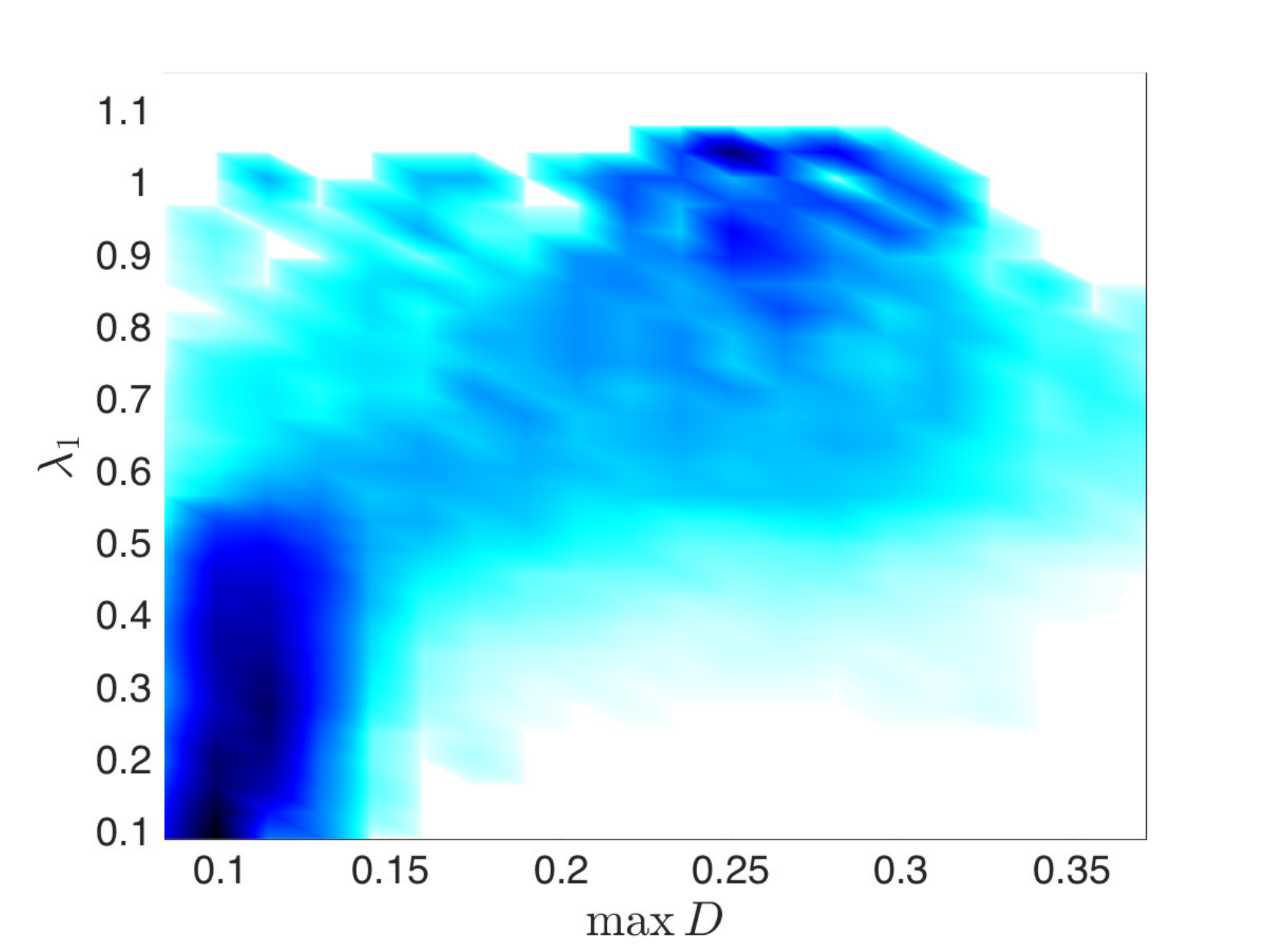}}
\subfigure[]{\includegraphics[width=.32\textwidth]{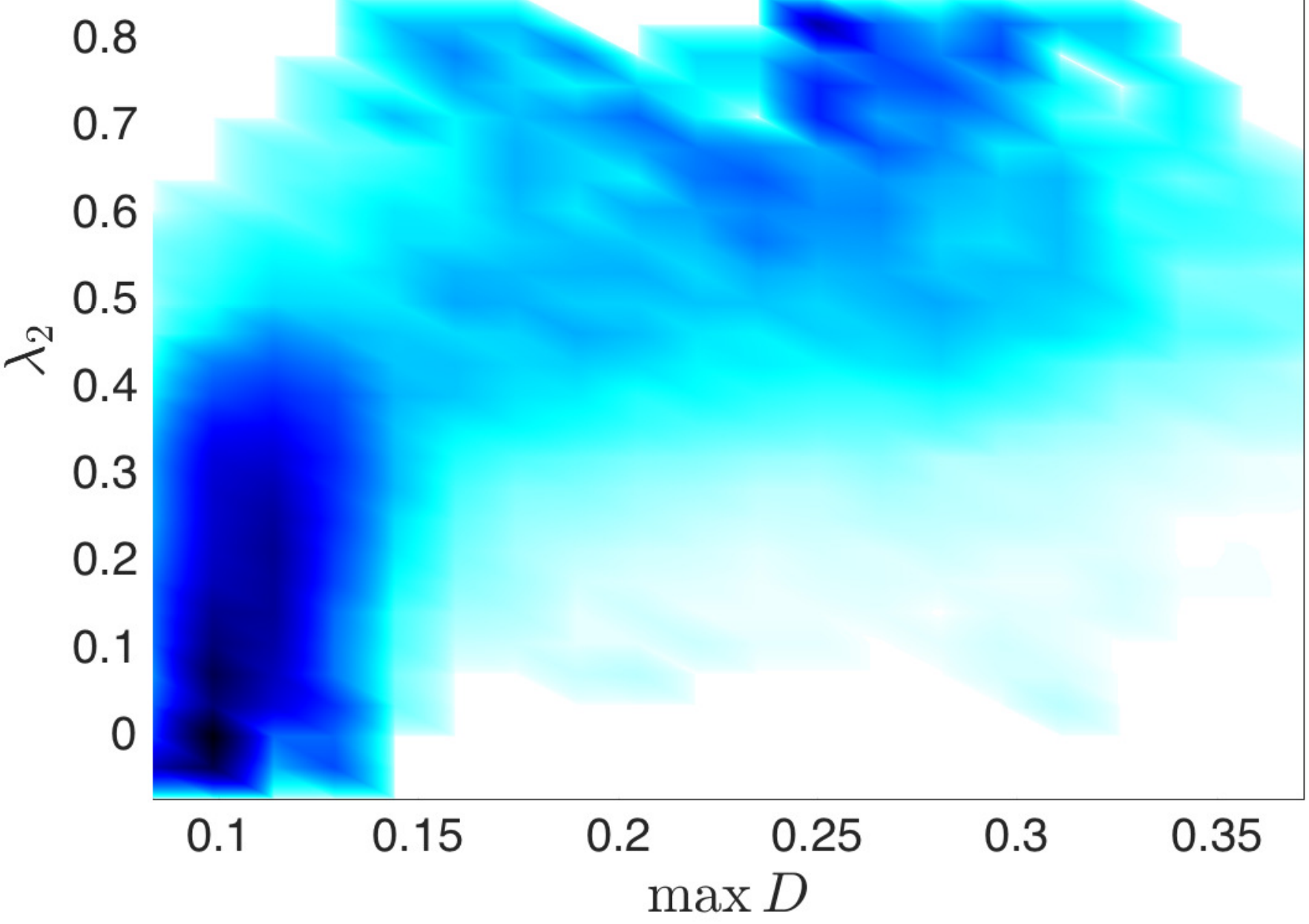}}
\subfigure[]{\includegraphics[width=.32\textwidth]{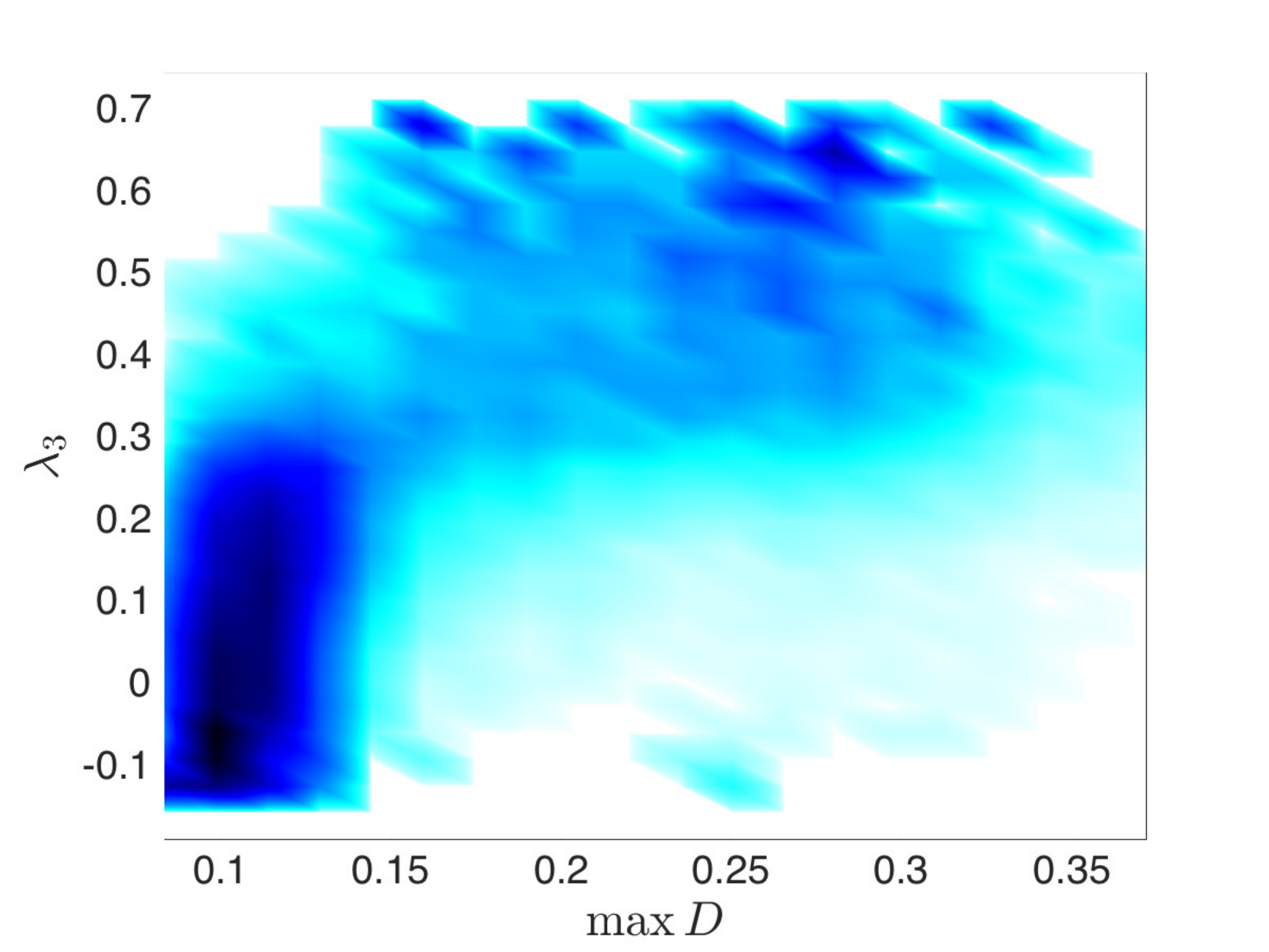}}
\caption{Conditional PDF of the first three eigenvalues of $\vc S_8$ 
and the maximal dissipation $\max_\tau D$ where the maximum is taken
over $\tau\in[t+t_i,t+t_f]$ with $t_i=3$ and $t_f=4$.} 
\label{fig:cond_pdf_lam123}
\end{figure}

To this end, we use the joint probability density function (PDF) of $\bar q$ and $\alpha$. 
The joint PDF of $\bar q$ and $\alpha$ is defined as the scalar function $p_{\bar q,\alpha}:\mathbb R\times\mathbb R\to \mathbb R^+$ that satisfies
\begin{equation}
\mathcal P(q_1\leq \bar q \leq q_2, \alpha_1\leq \alpha\leq \alpha_2)=
\int_{q_1}^{q_2}\int_{\alpha_1}^{\alpha_2}p_{\bar q,\alpha}(\bar q',\alpha')\id\bar q'\,\id\alpha',
\label{eq:joint_pdf}
\end{equation}
for all $q_1,q_2, \alpha_1, \alpha_2\in \mathbb R$ where 
$\mathcal P$ denotes the probability. 
The conditional probability density function of $\bar q$ (conditioned on $\alpha$) is then given by
\begin{equation}
p(\bar q|\alpha)=\frac{p_{\bar q,\alpha}(\bar q,\alpha)}{p_{\alpha}(\alpha)},
\label{eq:cond_pdf}
\end{equation}
where $p_{\alpha}$ is the probability density function of the indicator $\alpha$.

Roughly speaking, $p(\bar q(t)=\bar q_0|\alpha(t)=\alpha_0)$ denotes the 
likelihood of the maximum of the scalar $q$ over the time interval $[t+t_i,t+t_f]$
being $q_0$ given that the value of $\alpha$ at time $t$ is $\alpha_0$.
More precisely, the conditional probability of $\bar q$ over the time interval $[t+t_i,t+t_f]$ 
being greater than a prescribed value $q_0$ is given by
\begin{align}
\mathcal P\Big(\bar q(t)>q_0|\alpha(t)=\alpha_0\Big) & =
\mathcal P\left( \max_{\tau\in[t+t_i,t+t_f]}q(\tau)\geq q_0| \alpha(t)=\alpha_0\right)\nonumber\\
& = \int_{q_0}^{\infty}p(\bar q'|\alpha_0)\id \bar q'.
\end{align}

In particular, if an extreme event corresponds to values of $q$ greater than a prescribed critical value $q_c$,
the probability of the extreme event taking place over the time interval $[t+t_i,t+t_f]$, given that $\alpha(t)=\alpha_0$,
is measured by
\begin{equation}
P_{EE}(\alpha_0):=\mathcal P\left( \max_{\tau\in[t+t_i,t+t_f]}q(\tau)\geq q_c| \alpha(t)=\alpha_0\right)
 = \int_{q_c}^{\infty}p(\bar q'|\alpha_0)\id \bar q',
\label{eq:Pee}
\end{equation}
where $P_{EE}$ denotes the probability of an extreme event taking place over the future time interval $[t+t_i,t+t_f]$.
\begin{figure}
\centering
\subfigure[]{\includegraphics[width=.3\textwidth]{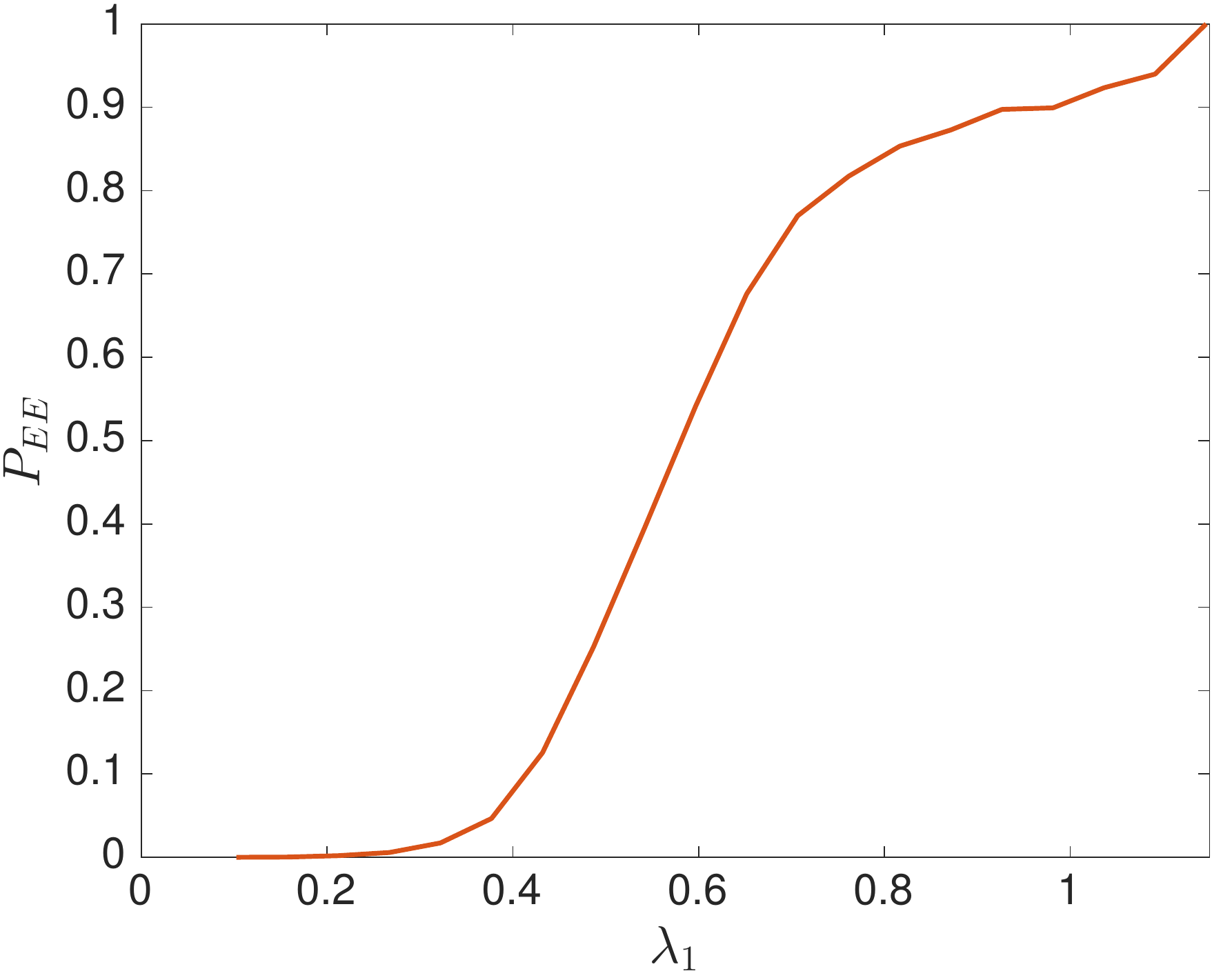}}
\subfigure[]{\includegraphics[width=.33\textwidth]{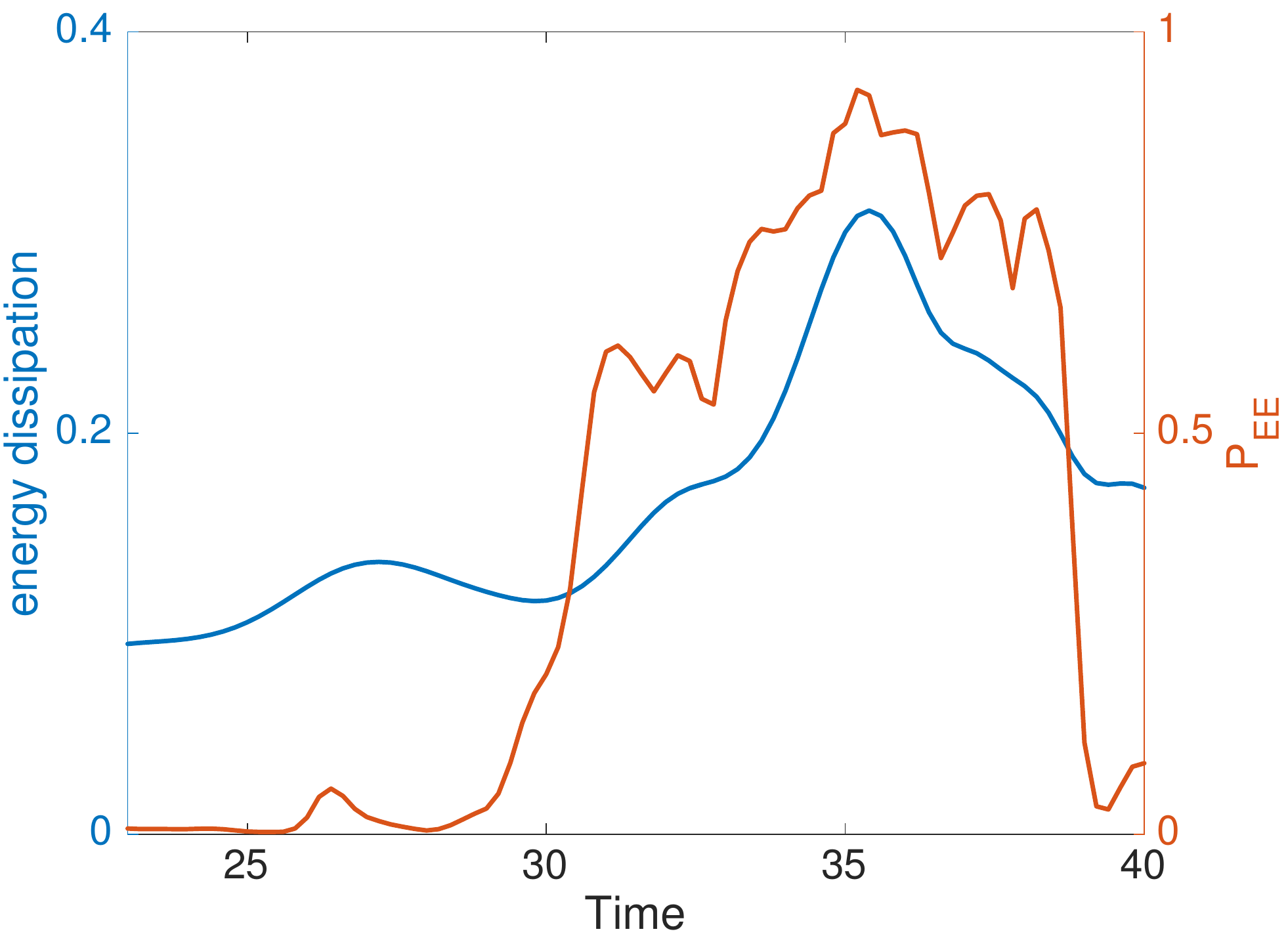}}
\subfigure[]{\includegraphics[width=.33\textwidth]{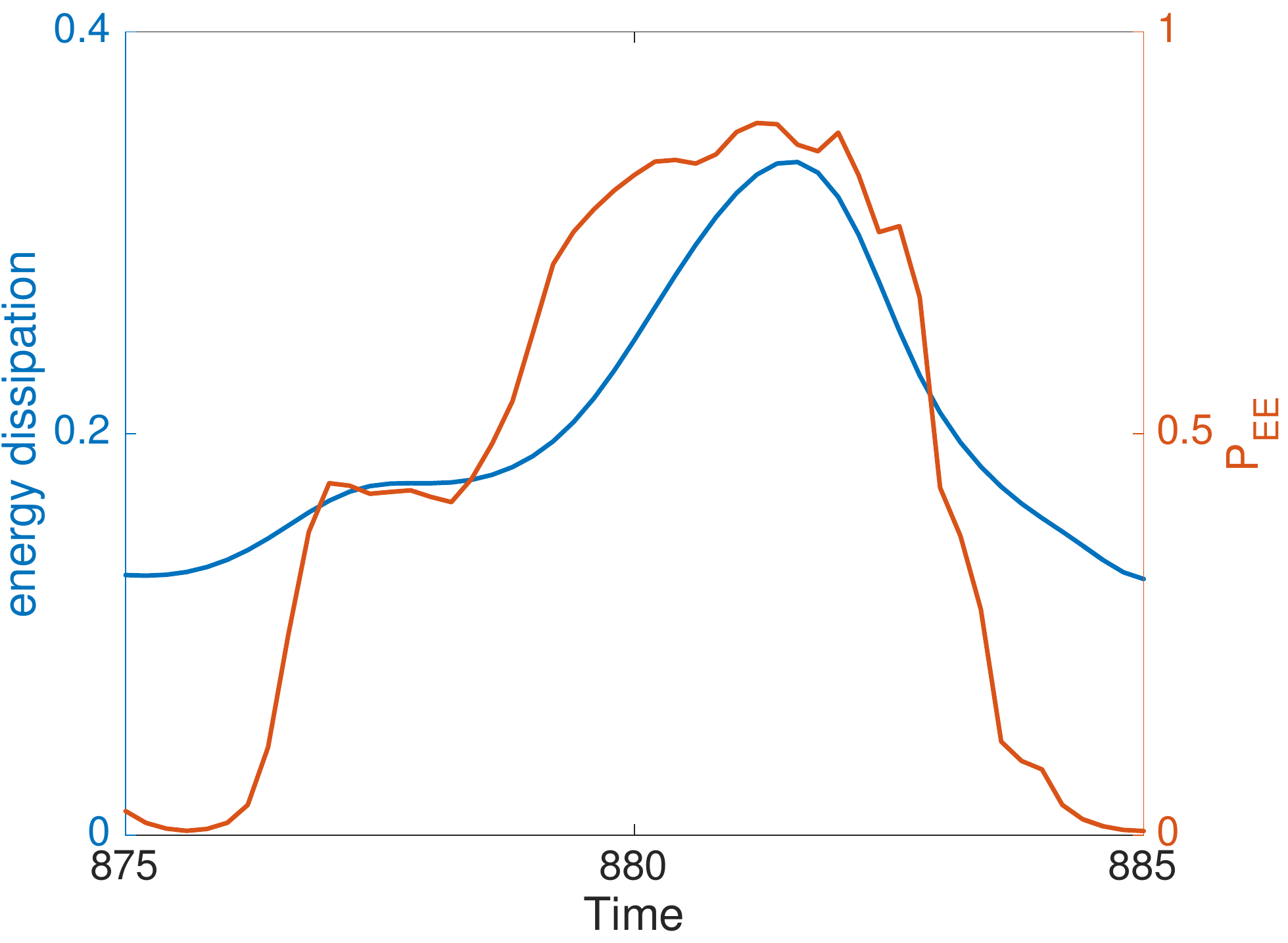}}
\caption{$t_i=3$ and $t_f=5$. Left: Probability of the extreme energy dissipation $P_{EE}$ as a function 
of the value of the indicator $\lambda_1$.
} 
\label{fig:Pee}
\end{figure}

In the case of the Kolmogorov flow, the observed quantity $q$ is the energy dissipation $D$ and
the indicator $\alpha$ is one the eigenvalues $\lambda_{i}$ of the reduced
symmetric tensor $\vc S_r$ (see equation~\eqref{eq:Sr}) with $r=8$. The joint PDF $p_{\bar D,\lambda_i}$
and the PDF $p_{\lambda_i}$ are approximated from a large set of 
numerical simulations containing roughly $85,000$ data points. 
The conditional PDF $p(\bar D|\lambda_i)$ then is computed through
the Bayesian relation~\eqref{eq:cond_pdf}.

Figure~\ref{fig:cond_pdf_lam123} shows the resulting conditional PDF $p(\bar D|\lambda_i)$
for the three largest eigenvalues of $\vc S_r$. As the three conditional PDFs are qualitatively 
similar, we will only discuss the one corresponding to the largest eigenvalue $\lambda_1$.

The conditional PDF exhibits a `bimodal' structure. For $0<\lambda_1(t)< 0.55$, 
the maximal future value of the energy dissipation $\max_{\tau \in[t+t_i,t+t_f]}D(\tau)$ 
is most likely to lie between $0$ and $0.15$ (the lower left dark region in figure~\ref{fig:cond_pdf_lam123}(a)).
A sharp transition is observed for $0.55<\lambda_1(t)$ such that for this range of the eigenvalue $\lambda_1$,
the energy dissipation is more likely to assume values larger that $0.15$ over the future time interval
$[t+t_i,t+t_f]$.

Using this conditional PDF, we compute the probability of extreme events $P_{EE}$ from equation~\eqref{eq:Pee}. 
From the time series presented in figure~\ref{fig:D_1001}, 
it is reasonable to associate a burst with values of the energy dissipation larger than $0.2$. 
We use this value as the critical energy dissipation (i.e. $D_c=0.2$) above which 
an extreme event is recorded. The resulting probability function is plotted in figure~\ref{fig:Pee}(a).
If at a time instant $t$, the value of $\lambda_1$ is smaller than $0.4$, the probability of 
$D(\tau)>D_c$ over the future time interval $\tau\in[t+t_i,t+t_f]$ is virtually zero. For larger values of 
$\lambda_1$, the probability of an extreme event increases monotonically. At $\lambda_1=0.55$,
the probability of an upcoming extreme event is greater than $50\%$.
Eventually, this probability grows to above $80\%$ at $\lambda_1\simeq 0.8$.

Using the computed probability of extreme event $P_{EE}$, we predict, at 
every given time $t$, the probability that an extreme event takes place over the future time interval $[t+t_i,t+t_f]$.
Figure~\ref{fig:Pee} (panels (b) and (c)), shows two select time windows 
over which an extreme event occurs. Away from the extreme event, the probability $P_{EE}$
is very low. Just before the extreme event, this probability grows predicting the
upcoming extreme events at least $t_i=3$ time units in advance. 
\begin{figure}
\centering
\subfigure[]{\includegraphics[width=.4\textwidth]{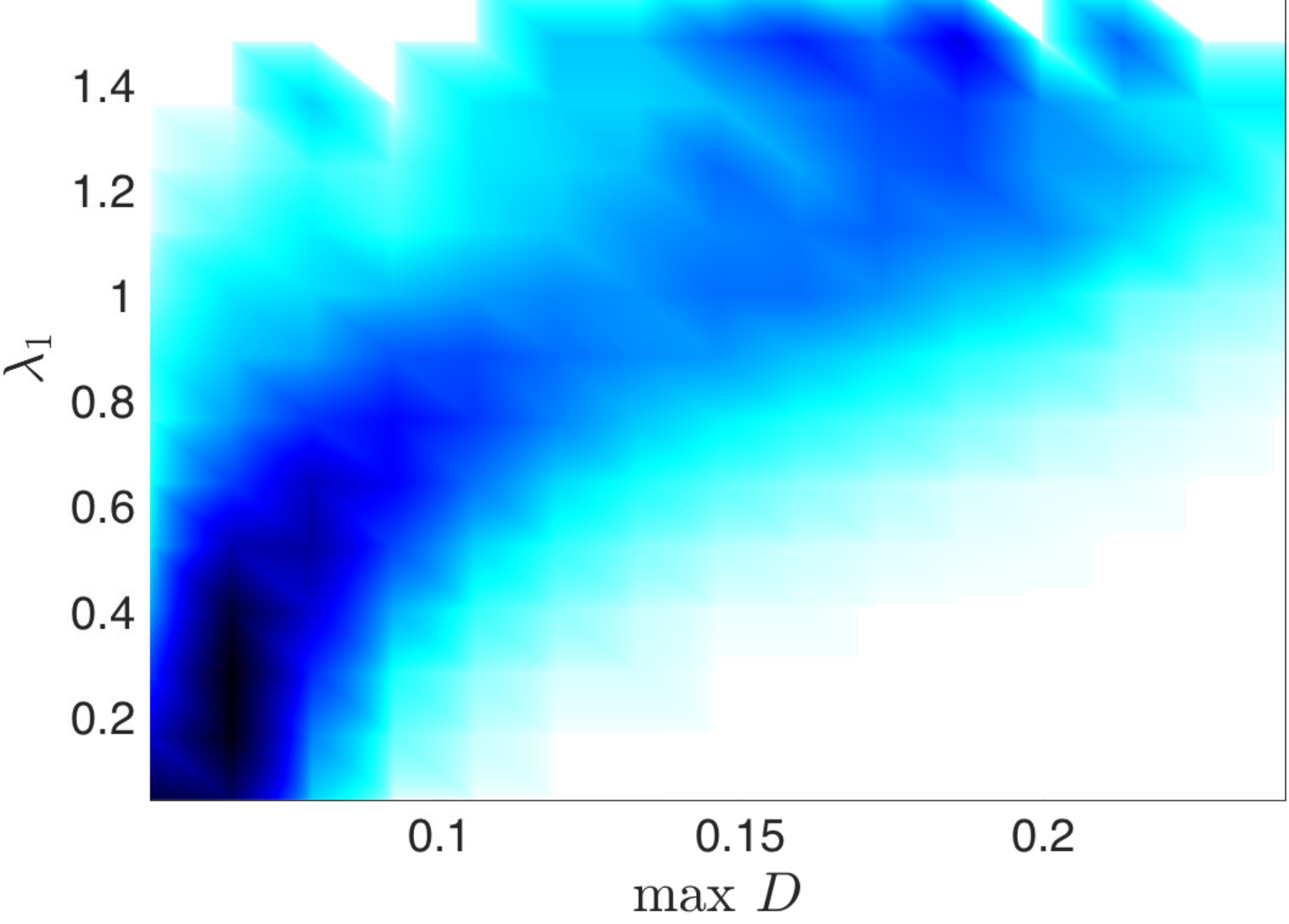}}
\subfigure[]{\includegraphics[width=.4\textwidth]{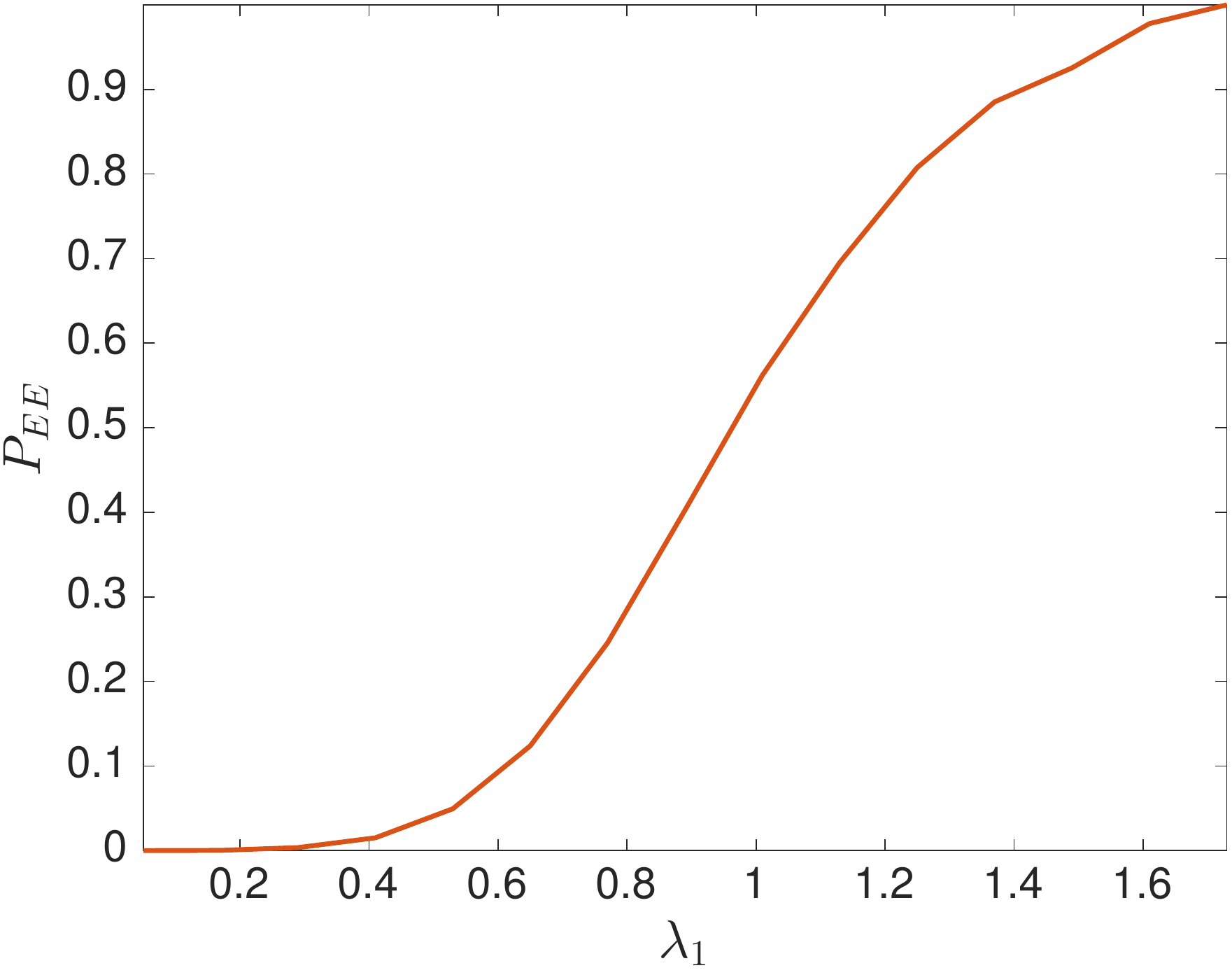}}
\caption{Conditional PDF (a) and the probability of upcoming extreme energy
dissipation (b) for Reynolds number $Re=100$.} 
\label{fig:R100}
\end{figure}

While the above results are reported at $Re=40$, we point out that similar
conclusions hold at higher Reynolds numbers. Figure~\ref{fig:R100}, for instance,
shows the conditional PDF and the probability of extreme events at $Re=100$.
To fully resolve the flow, the higher resolution of $256\times 256$ Fourier modes 
are used at this $Re$ number. On the other hand, to 
keep the computational cost reasonable, the linearized operator
is reduced to four OTD modes, i.e., $r=4$.

%
\subsection{Comparison with dynamic mode decomposition}
We carry out a caparison in this section to highlight that 
the correct choice of the modes to which the linear operator $\lins{\vc u}$
is reduced is essential. To this end, we repeat the analysis of Section~\ref{sec:pdf},
but this time we reduce the operator $\lins{\vc u}$ to the modes 
obtain from Dynamic Mode Decomposition (DMD).
DMD was proposed by~\citet{schmid10} for extracting
a linear approximation to the flow map of a nonlinear dynamical system.
The resulting \emph{dynamic modes} (or DMD modes) have proven insightful in the analysis of
fluid flows~\cite{schmid11,schmid11no2} and shown to have intricate connections to the Koopman and Fourier 
modes of time periodic solutions~\cite{rowley09,chen12}.
\begin{figure}
\centering
\includegraphics[width=.5\textwidth]{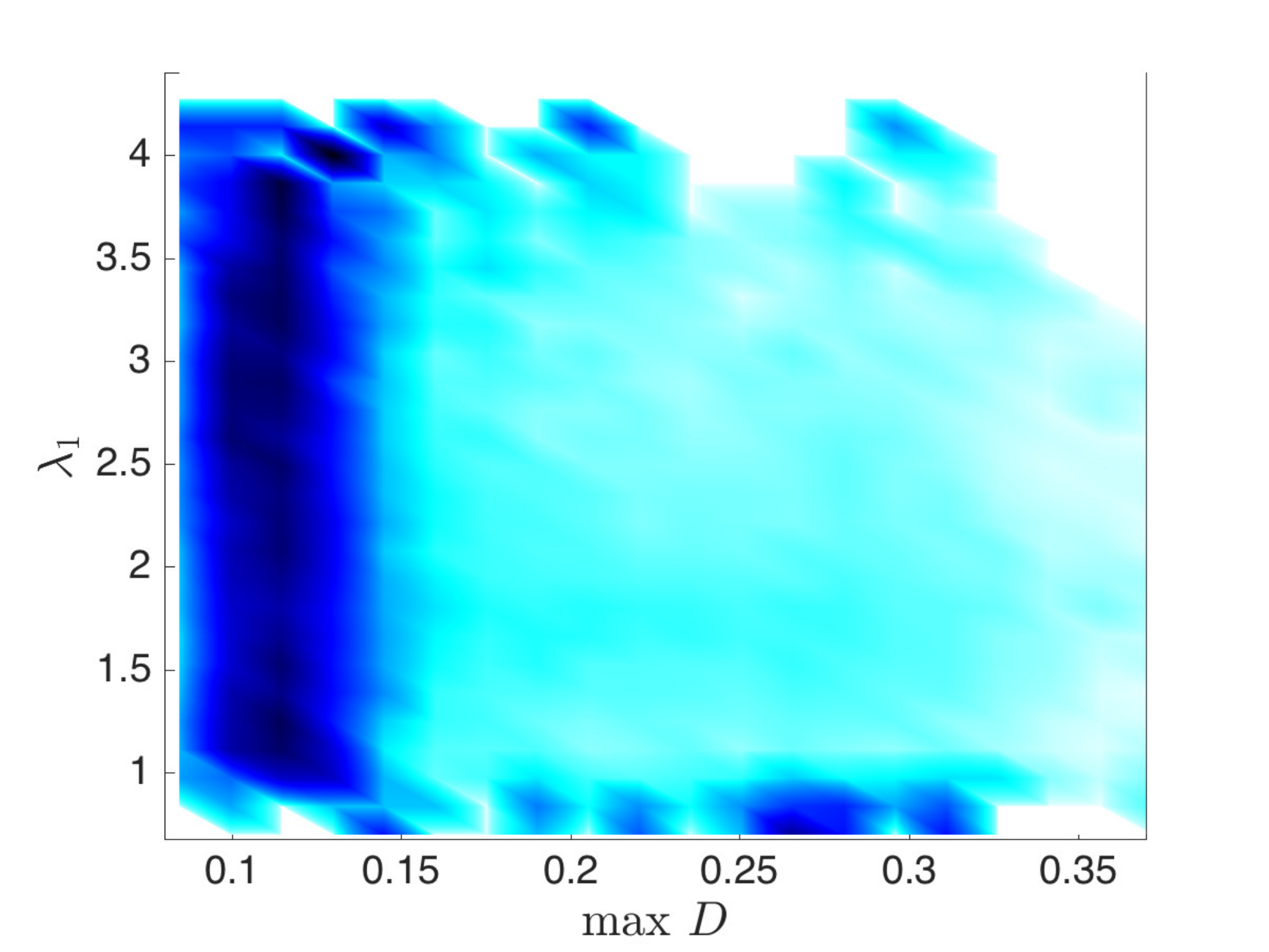}
\caption{Same as figure~\ref{fig:cond_pdf_lam123}(a) but now the linear 
operator is reduced to the eight most dominant DMD modes.
}
\label{fig:pdf_dmd}
\end{figure}

Since the DMD modes are not flow invariant (see Definition~\ref{def:flowInv}),
the reduction of the linear operator $\lins{\vc u}$ to these modes is not dynamically meaningful.
As a result, the eigenvalues of the symmetric tensor reduced to DMD modes are not expected to
reflect the true growth (or decay) of perturbations.
To illustrate this, we use the algorithm introduced by~\citet{schmid10} to compute DMD modes
from $500$ sequential snapshots of the Kolmogorov flow, each $0.2$ time units apart.
Next we restrict the operator $\lins{\vc u}$ to the eight most dominant DMD modes
and compute the largest eigenvalue of its symmetric part along all previously computed turbulent 
trajectories $\vc u(t)$. The resulting conditional PDF is shown in figure~\ref{fig:pdf_dmd}. 
As opposed to the OTD modes (cf. figure~\ref{fig:cond_pdf_lam123}),
the extreme episodes of the energy dissipation do not show a signature 
in the DMD-reduced operator.

\section{Spatially localized extreme events}\label{sec:mnls}
The energy dissipation in turbulent flows, as discussed in Section~\ref{sec:kolm}, is a global
feature of the state. In spatiotemporal chaos, however, local rare extreme events, 
in the form of spatially localized structures, are possible . A famous example of 
such localized extreme events is the ocean rogue waves. 
Such localized phenomena cannot be quantified from global quantities such as 
the eigenvalues of the linear operator.

In this section, we illustrate that localized features of the OTD modes can still
be of significance for the analysis of spatially localized extreme events. 
To illustrate this, we consider the modified nonlinear Schr\"odinger (MNLS) equation which
is an approximation to the evolution of sea surface elevation in deep waters~\citep{dysthe79}.
The MNLS equation is a higher order perturbative approximation compared to the 
nonlinear Schr\"odinger equation derived by~\citet{zakharov68}.
Recently, more quantitative methods for the analysis of the extreme waves in the MNLS equation
have been developed~\citep{cousins15,cousins16,mohamad15}.

\subsection{MNLS equation}
For a complex valued function $u(x,t)$, the MNLS equation (in dimensionless variables) reads
\begin{equation}
\partial_t u =F(u),
\end{equation}
with
\begin{equation}
F(u)=-\frac{1}{2}\partial_x u -\frac{\i}{8}\partial_x^2 u+\frac{1}{16}\partial_x^3 u-\frac{\i}{2}|u|^2u
-\frac{3}{2}|u|^2\partial_x u-\frac{1}{4}u^2\partial_xu^\ast-\i\, u\Phi(u),
\label{eq:mnls}
\end{equation}
where $\i=\sqrt{-1}$, $x\in[0,L]$ and $u(x,t)\in \mathbb C$. The asterisk sign denotes the complex
conjugation. The function $\Phi$ is derived from the velocity potential 
$\phi$,
\begin{equation}
\Phi(u):=\partial_x\phi\Big|_{z=0}=-\frac{1}{2}\mathcal F^{-1}\left[|k|\mathcal F[|u|^2]\right],
\label{eq:mnls_Phi}
\end{equation}
where $\mathcal F$ denotes the Fourier transform.
The modulus $|u(x,t)|$ is the wave envelope for the surface elevation $h(x,t)$. To the leading 
order approximation, we have $h(x,t)=\mbox{Re}\left[u(x,t)\exp(\i (x-t)) \right]$.
\begin{figure}
\centering
\includegraphics[width=\textwidth]{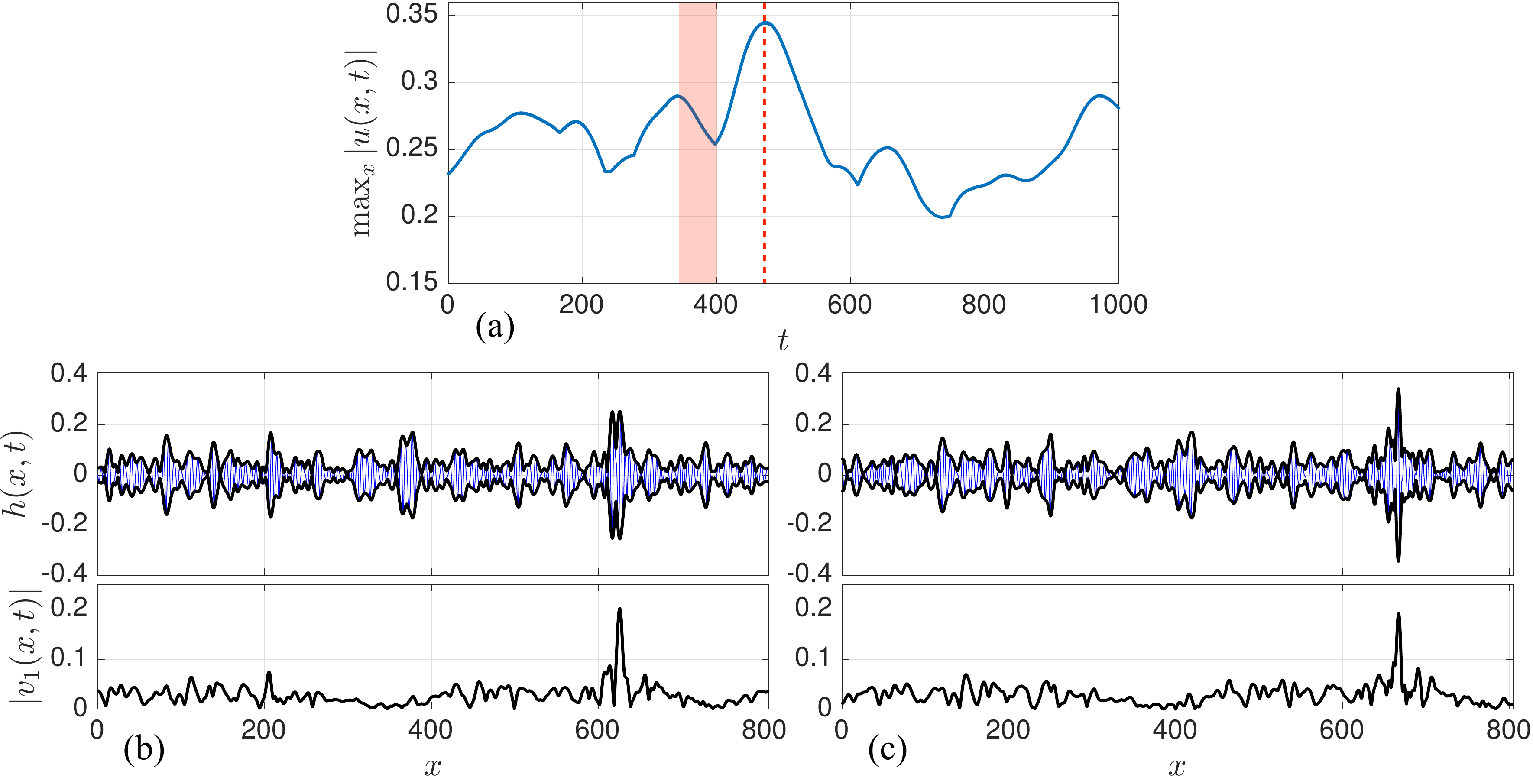}
\caption{
(a) The spatial maximum of $|u|$ as a function of time $t$. An extereme event occurs at around
$t=475$ where $\max_x|u|\simeq 0.34$.
(b, c) The surface elevation $h(x,t)$ (blue color) and and the modulus of the OTD mode $|v_1|$ 
at times $t=400$ (b) and $t=475$ (c). The thick black curves in the plots of $h(x,t)$ mark the envelopes
$\pm |u(x,t)|$.
}
\label{fig:mnls_wave}
\end{figure}

We solve the MNLS equation with the 
initial conditions $u(x,0)=u_0(x)$ with Gaussian energy spectra and random phases. 
More precisely,
the Fourier transform of the initial condition is given by
\begin{equation}
\widehat{u_0}(k)=\sqrt{2\frac{2\pi}{L}N(q_k)} e^{\i\theta_k},
\label{eq:ic_mnls}
\end{equation} 
where 
\begin{equation}
N(q_k):=\frac{\epsilon^2}{\sigma\sqrt{2\pi}}e^{-\frac{q_k^2}{2\sigma^2}},
\end{equation}
is a normal distribution, $\theta_k$ are random phases uniformly distributed over $[0,2\pi]$
and $q_k=2\pi k/L$ is the wave number over the periodic domain of length $L$. There are three free parameters:
$\epsilon$ that controls the wave height,
$\sigma$ which is the standard deviation of the Gaussian distribution and controls the width of
the spectrum of the wave and finally $L$ which is the length of the periodic domain, $x\in[0,L]$.

It is well-known that the Gaussian wave groups~\eqref{eq:ic_mnls}
can grow due to the Benjamin-Feir instability~\citep{benjamin67}
to form extreme waves. The Benjamin-Feir Index (BFI)
$2\sqrt{2}\epsilon/\sigma$ provides an indicator for the probability of the extreme waves 
taking place. For large enough BFI, the nonlinear terms dominate, leading to large amplitude 
waves~\cite{janssen03}. If BFI is too large, however, the extreme waves happen quite often. 
To realize rare extreme waves, therefore, a moderate BFI value should be used.
Following~\citet{mohamad15}, we use the parameter values $\epsilon=0.05$,
$\sigma=0.2$ and $L=256\pi$, resulting in BFI$=0.71$. This BFI value 
allows for the formation of extreme waves at a moderate frequency (not too often
and not too rare).

We solve the MNLS~\eqref{eq:mnls} equation and its associated OTD equation~\eqref{eq:otd_pde} where 
$\langle\cdot,\cdot\rangle$ now denotes the standard $L^2$ inner product on 
complex valued functions,
\begin{equation}
\langle v,w\rangle :=\int_{0}^Lv(x)w^\ast(x)\id x.
\end{equation}
The initial condition for the OTD modes are sinusoidal and are given by
$$v_i(x,0)=\sqrt{\frac{2}{L}}\sin\left(\frac{2\pi i}{L}x\right).$$
The computation of the OTD modes requires the linearization of 
the operator~\eqref{eq:mnls} as outlined in Appendix~\ref{app:linF}.

For the numerical integration of the MNLS equation (and its associated OTD equation),
we use a second-order exponential time differencing scheme~\cite{cox02,berland07}
in time and a pseudo-spectral scheme for evaluating the spatial derivatives with $2^{11}$ Fourier modes.
For the statistical analysis presented in the next section, we compute $200$
MNLS trajectories, each of length $1000$ time units,
from the initial conditions of the form~\eqref{eq:ic_mnls}.

\subsection{Extreme waves and the OTD modes}
Figure~\ref{fig:mnls_wave} shows a time window over which an extreme wave appears at 
around $t=475$ with a wave height of approximately $0.34$ (see panel (a)).
Panel (b) shows a snapshot of the wave, $75$ time units earlier at $t=400$.
It exhibits a twin wave packet at around $x=610$. Whether this twin 
wave packets lead to an extreme wave depends on the energies and the phases
of the packets. A simple extrapolation will rule out the possibility of an extreme wave 
since the wave height has been decaying over the last $50$ time units (the red shaded
area in figure~\ref{fig:mnls_wave}(a)).
\begin{figure}
\centering
(a)\includegraphics[width=.45\textwidth]{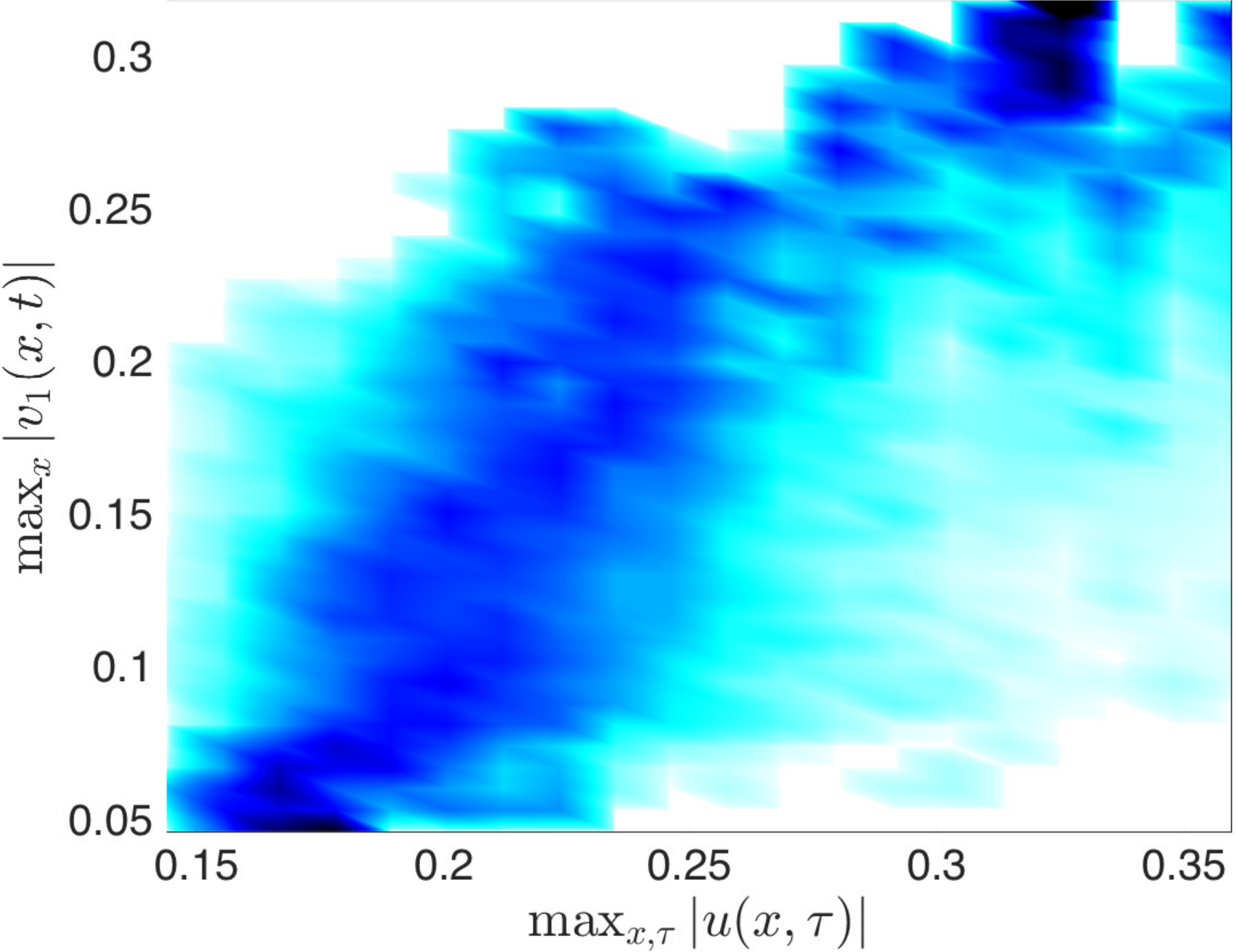}
(b)\includegraphics[width=.45\textwidth]{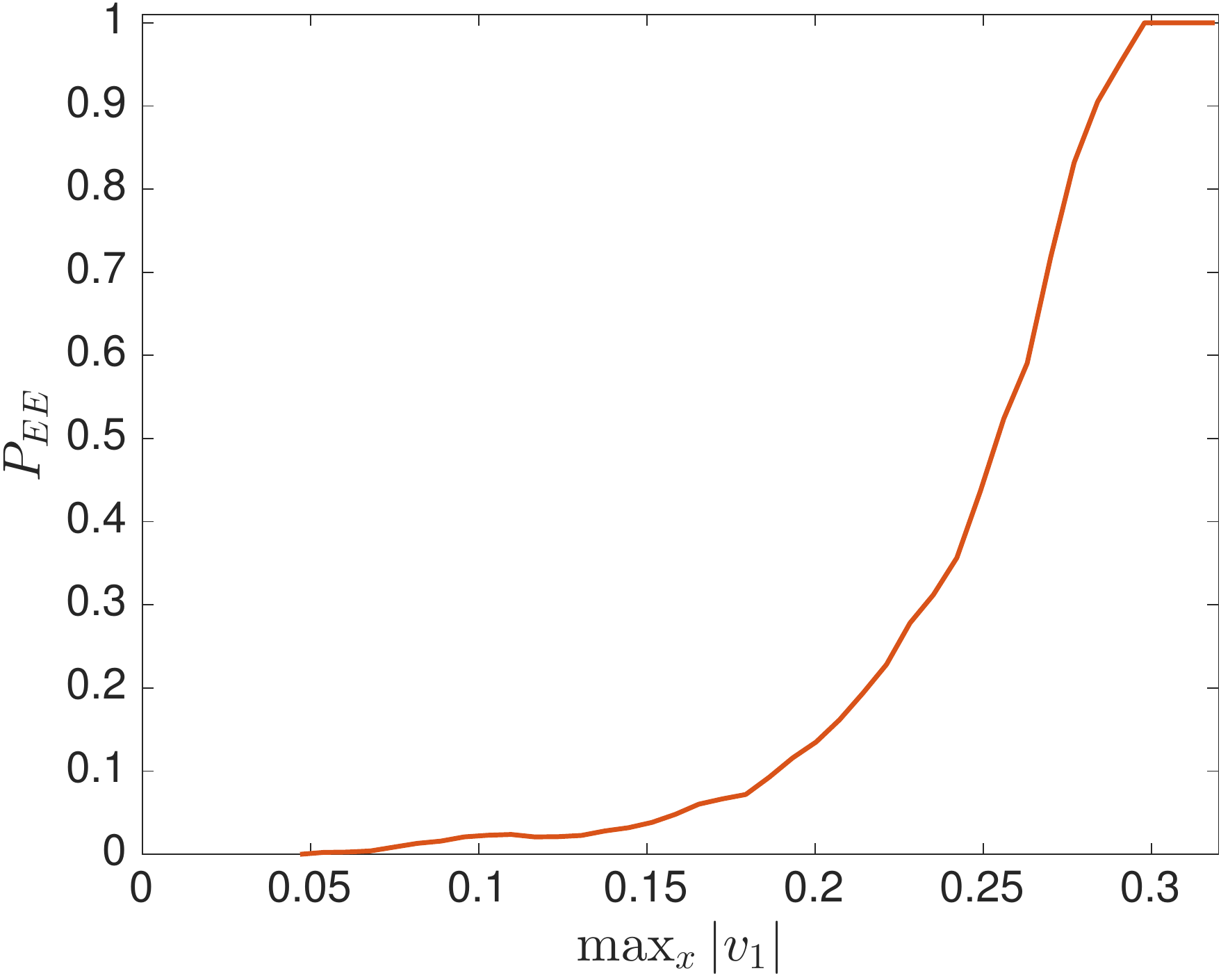}
\caption{
(a) Conditional PDF for the maximum modulus of the OTD mode $v_1$
and the solution of the MNLS equation. The maxima are taken over $x\in[0,L]$
and $\tau\in[t+t_i,t+t_f]$ with $t_i=25$ and $t_f=26$.
(b) The probability of an extreme event $P_{EE}$ computed from the conditional PDF.
}
\label{fig:mnls_pdf}
\end{figure}

During this decay period, however, the OTD mode $v_1$ shows a persistent 
localized peak at the same location as the twin wave packets. This signals
a persistent localized instability that grows to lead to the extreme wave
at time $t=475$ as shown in figure~\ref{fig:mnls_wave}(c).

As in the case of the Kolmogorov flow, we use Bayesian statistics 
to quantify the relation between extreme MNLS waves and the 
localized peaks of the associated OTD modes. Based on the foregoing observation,
we use the maximum height of the 
first OTD mode $v_1$ as the indicator $\alpha$. The quantity to be predicted is the
maximum height of the MNLS solution $u$. More precisely,
$$q(t)=\max_{x\in[0,L]}|u(x,t)|,\quad \alpha(t)=\max_{x\in[0,L]}|v_1(x,t)|.$$
The conditional PDF $p(\bar q|\alpha)$ is computed as in Section~\ref{sec:pdf}.
For a given critical wave height $h_c$, 
the probability of the rare event is given as in equation~\eqref{eq:Pee}
by 
\begin{equation}
P_{EE}(\alpha_0):=
\mathcal P\left( \max_{\tau\in[t+t_i,t+t_f]}\max_{x\in[0,L]}|u(x,\tau)|\geq h_c\vline
\max_{x\in[0,L]}|v_1(x,t)|=\alpha_0\right).
\end{equation}

Figure~\ref{fig:mnls_pdf} shows the conditional PDF $p(\bar q|\alpha)$ and
the probability of an extreme wave with the critical wave height $h_c=0.28$.
This critical wave height is approximately the mean plus two standard deviation
of $\max_x|u|$ for all the data computed. 

\section{Summary and conclusions}\label{sec:concl}
We proposed operational indicators for the prediction of 
rare extreme events (or bursts) in high dimensional dynamical systems. 
The motivation for our indicators is based on the observations
made about slow-fast systems where the bursts 
occur along orbits that are transverse and homoclinic to
the slow manifold~~\citep{haller99,guck12,guck15}. This geometric picture does not lead to an operational 
method in complex high-dimensional systems where
a clear separation between the slow and fast variables is unavailable~\cite{pope06}.

We showed that for such systems a signature of bursting can be traced in the
eigenvalues of the symmetric part of the linearized dynamics. More precisely, 
we use the largest eigenvalue $\lambda_1$ of the symmetric part of the 
linearized operator as our indicator. Computing these
eigenvalues in high dimensional systems is computationally expensive. 
Thanks to the recently
introduced notion of Optimally Time Dependent (OTD) modes~\cite{otd}, however,
one can reduce the linear operator, in a dynamically consistent fashion, to
its most unstable subspace. The reduced operator is low dimensional 
and its invariants can be readily computed. 

We devised a low dimensional ODE in Section~\ref{sec:bursting} which has an unambiguous 
bursting mechanism. For this simple model we showed that the eigenvalue $\lambda_1$
becomes uniformly positive several time units before the burst. This allows for 
instantaneous perturbations within the corresponding subspace to grow.
Moreover, the OTD mode aligns with the direction of the growth (i.e. orthogonal to the $x-y$ plane). 
These together successfully predict the upcoming extreme event.

In the body forced Navier--Stokes equation considered in Section~\ref{sec:kolm}, 
the situation is more complicated as the symmetric part of the reduced operator has several 
eigenvalues that are positive for all times. The largest eigenvalue $\lambda_1$, however,
increases significantly before a burst in the energy dissipation takes place. 
Using Bayesian statistic, we showed that large values of the eigenvalue $\lambda_1$
do in fact predict upcoming bursts in the energy dissipation. While the results are presented for prediction time
$t_i=3$, they are robust to small variations of this time window. 
If the prediction time is set too large (larger than $t_i=5$, here), however, the indicator fails to predict the  
bursts. The predictability time, of course, is problem dependent and is expected to be inversely proportional
to the dominant Lyapunov exponent of the system~\cite{broer2010}.

We also considered extreme waves in a unidirectional model of the 
nonlinear surface waves in deep ocean. 
As opposed to the energy dissipation in Navier--Stokes equations, 
extreme waves are localized
in space. Therefore, we do not expect the eigenvalue $\lambda_1$
(as a global quantity) to bear significance in their creation.
We observe instead that the most unstable OTD mode localizes and grows 
before an extreme wave appears. The spatial location where the OTD mode
localizes is precisely where the extreme wave occurs later in time. 
This observation indicates a promising direction for space-time prediction
of the extreme water waves, complementing the recent work of~\citet{cousins15,cousins16}.

Finally, we point out that the OTD modes are instrumental to the evaluation of our indicators.
This imposes an additional computational cost as the OTD equations need to be solved simultaneously
with the governing equations. Moreover, it necessitates that a model of the system 
is available as a set of differential equations. Therefore, modifying the indicator so that it is applicable to
model-independent predictions is highly desirable.
Future work also involves the application of the presented ideas on the filtering and prediction of stochastic dynamical systems exhibiting rare events. 
\subsection*{Acknowledgments}
The authors have been supported through the Air Force Office of Scientific Research (AFOSR YIP 16RT0548), the Office of Naval Research (ONR YIP N00014-15-1-2381), and the National Science Foundation (NSF EAGER ECCS 15-1462254).
\begin{appendices}
\section{The linearization of the MNLS equation}\label{app:linF}
We denote the linearization (or G\^ateaux differential) 
of the differential operator $F$ defined in equation~\eqref{eq:mnls} by 
$L_{u}(\cdot)$ which reads
\begin{align}
L_{u}(v) :=&\lim_{\epsilon\to 0}\frac{F(u+\epsilon v)-F(u)}{\epsilon}\nonumber\\
 =& -\frac{1}{2}\partial_x v -\frac{\i}{8}\partial_x^2 v+\frac{1}{16}\partial_x^3 v\nonumber\\
   & -\frac{\i}{2}\big(2|u|^2v+u^2v^\ast\big) \nonumber\\
   & -\frac{3}{2}\big(u^\ast v \partial_x u + uv^\ast \partial_x u+|u|^2\partial_x v\big)\nonumber\\
   & -\frac{1}{4}\big(2uv\partial_x u^\ast+u^2\partial_x v^\ast\big)\nonumber\\
   &+\frac{\i}{2}u\,\mathcal F^{-1}\big[|k|\,\mathcal F[uv^\ast+vu^\ast]\big]+\frac{\i}{2}v\,\mathcal F^{-1}\big[|k|\mathcal F(|u|^2)\big].
\label{eq:mnls_L}
\end{align}
The only nontrivial calculation above is the last line, corresponding to the linearization 
of the term $u\Phi(u)$ in~\eqref{eq:mnls}, which we detail below. First we note that
\begin{equation}
(u+\epsilon v)\Phi(u+\epsilon v)-u\Phi(u) = \epsilon u\, \id\Phi(u;v)+\epsilon v \Phi(u)+\mathcal O(\epsilon^2),
\end{equation}
where 
\begin{equation}
\id\Phi(u;v)=\lim_{\epsilon\to 0}\frac{\Phi(u+\epsilon v)-\Phi(u)}{\epsilon}.
\end{equation}
From the definition of $\Phi$ (see equation~\eqref{eq:mnls_Phi}), we have
\begin{align}
\mathcal F[\Phi(u+\epsilon v)] & =-\frac{1}{2}|k|\,\mathcal F[|u+\epsilon v|^2]\nonumber\\
&=\mathcal F[\Phi(u)]-\epsilon\frac{1}{2}|k|\, \mathcal F[uv^\ast+vu^\ast]+\mathcal O(\epsilon^2),
\end{align}
which yields
\begin{equation}
\id\Phi(u;v)=-\frac{1}{2}\mathcal F^{-1}\big[|k|\,\mathcal F[uv^\ast+vu^\ast]\big].
\end{equation}
This completes the derivation of~\eqref{eq:mnls_L}.
\end{appendices}

%
\end{document}